\documentclass[10pt]{amsart}
\usepackage{amssymb}

\newtheorem{lemma}{Lemma}[section]

\newtheorem{theorem}[lemma]{Theorem}
\newtheorem{remark}[lemma]{Remark}
\newtheorem{proposition}[lemma]{Proposition}
\newtheorem{corollary}[lemma]{Corollary}
\newtheorem{definition}[lemma]{Definition}

\begin{document}

\title{On the topology of non-isolated real singularities}

\author{Nicolas Dutertre}
\address{Laboratoire angevin de recherche en math\'ematiques, LAREMA, UMR6093, CNRS, UNIV. Angers, SFR MathStic, 2 Bd Lavoisier 49045 Angers Cedex 01, France.}
\email{nicolas.dutertre@univ-angers.fr}
 
\thanks{The author is partially supported by the ANR project LISA 17-CE400023-01 and by the Funda\c{c}\~ao de Amparo \`a Pesquisa do Estado de S\~ao Paulo - FAPESP, Brazil}

\subjclass[2010]{32B05, 58K05, 58K65}

\keywords{Topological degree, Euler characteristic, Real Milnor fibres}

\begin{abstract}
Khimshiashvili proved a topological degree formula for the Euler characteristic of the  Milnor fibres of a real function-germ with an isolated singularity. We give two generalizations of this result for non-isolated singularities. As corollaries we obtain an algebraic formula for the Euler characteristic of the fibres of a real weighted-homogeneous polynomial and a real version of the L\^e-Iomdine formula. We have also included some results of the same flavor on the local topology of locally closed definable sets.
\end{abstract}

\maketitle
\markboth{N. Dutertre}{On the topology of non-isolated real singularities }

\section{Introduction}
Let $f:(\mathbb{R}^n, 0) \to (\mathbb{R},0)$ be an analytic function-germ with an isolated critical point at the origin. Khimshiashvili \cite{Khimshiashvili} proved the following formula for the Euler characteristic of the real Milnor fibres of $f$:
$$\chi \big( f^{-1}(\delta) \cap B_\epsilon \big) = 1-{\rm sign} (-\delta)^n {\rm deg}_0 \nabla f,$$
where $0 < \vert \delta \vert \ll \epsilon \ll 1$, $B_\epsilon$ is the closed ball centered at the origin of radius $\epsilon$ and ${\rm deg}_0 \nabla f$ is the topological degree of the mapping $\frac{\nabla f}{\vert \nabla f \vert} : S_\epsilon \to S^{n-1}$ (here $S_\epsilon$ is the boundary of $B_\epsilon$). Later Fukui \cite{Fukui} generalized this result for the fibres of a one-parameter deformation of $f$. A corollary of the Khimshiashvili formula due to Arnol'd \cite{Arnold} and Wall \cite{Wall} states that 
$$\chi(\{f \le 0 \} \cap S_\varepsilon)= 1- \hbox{\rm deg}_0 \nabla f, $$  $$\chi(\{f \ge 0 \} \cap S_\varepsilon)= 1+(-1)^{n-1}\hbox{\rm deg}_0 \nabla f,$$
and if $n$ is even,  
$$\chi(\{f = 0 \} \cap S_\varepsilon )= 2-2\ \hbox{\rm deg}_0 \nabla f.$$
In \cite{SzafraniecTopology86} Szafraniec extended the results of Arnold and Wall to the case of an analytic function-germ $f :(\mathbb{R}^n,0) \to (\mathbb{R},0)$ with non-isolated singularities. Namely he constructed two function-germs $g_-$ and $g_+$ with isolated critical points and proved that
$$\chi \left( \{f \le 0 \} \cap S_\epsilon \right) =1-{\rm deg}_0 \nabla g_+ \hbox{ and }  \chi \left( \{f \ge 0 \} \cap S_\epsilon \right) =1- {\rm deg}_0 \nabla g_-.$$
In \cite{SzafraniecGlasgow91} he improved this result for weighted homogeneous polynomials. If $f : \mathbb{R}^n \to \mathbb{R}$ is a weighted homogeneous polynomial then he constructed to polynomials $g_1$ and $g_2$ with an algebraically isolated critical point at $0$ such that 
$$\chi \left( \{f \le 0 \} \cap S^{n-1} \right) =1-{\rm deg}_0 \nabla g_1 \hbox{ and }  \chi \left( \{f \ge 0 \} \cap S^{n-1} \right) =1- {\rm deg}_0 \nabla g_2.$$
Thanks to the Eisenbud-Levine-Khimshiashvili formula \cite{EL,Khimshiashvili}, $\chi \left( \{f \le 0 \} \cap S^{n-1}\right)$ and $\chi \left( \{f \ge 0 \} \cap S^{n-1} \right)$ can be computed algebraically.

The aim of this paper is to extend the Khimshiashvili formula for function-germs with arbitrary singularities. We will work in the more general framework of definable functions. Let $f:(\mathbb{R}^n, 0) \to (\mathbb{R},0)$ be a definable function-germ of class $C^r$, $r \ge 2$. Our first new result is Lemma \ref{MilnorfibreLink} where we give a relation between the Euler characteristic of $f^{-1}(\delta) \cap B_\epsilon$ (resp. 
$f^{-1}(-\delta) \cap B_\epsilon$), with $0 < \delta \ll \epsilon \ll 1$, and the Euler characteristic of the link at the origin of $\{ f \le 0 \}$ (resp. $\{ f \ge 0 \}$). Applying the results of Szafraniec, we obtain our first generalization of the Khimshiashvili formula (Corollary \ref{KhimPolyBounded}) for polynomially bounded structures and an algebraic formula for the Euler characteristic of a regular fibre of a weighted homogeneous polynomial (Corollary \ref{FibreMilnorHomogeneous}). We note that the paper \cite{DimcaPaunescu} presents a different approach for the computation of this Euler characteristic.

Our second generalization of the Khimshiashvili formula is an adaptation to the real case of the methods based on the generic polar curve, introduced in the complex case by L\^e \cite{LeRelativePolar} and Teissier \cite{TeissierEqui,TeissierRelativePolar} and developed later by Massey \cite{MasseyInventiones1,MasseyInventiones2,MasseyLeCycles}. For $v \in S^{n-1}$, we denote by $\Gamma_v$ the following relative polar set:
$$\Gamma_v =\left\{ x \in \mathbb{R}^n \setminus \Sigma_f \ \vert \ {\rm rank}(\nabla f(x), v ) <2 \right\},$$
where $\Sigma_f =\left\{ x \in \mathbb{R}^n \setminus \Sigma_f \ \vert \  \nabla f(x) =0 \right\}$ is the critical locus of $f$.
For $v$ generic in $S^{n-1}$, $\Gamma_v$ is a curve. Let $\mathcal{B}$ be the set of its connected components. For each ${\rm \bf b} \in \mathcal{B}$, we denote by $\sigma({\rm \bf b})$ the sign of ${\rm det}\big[ \nabla f_{x_1},\ldots,\nabla f_{x_n} \big]$ on ${\rm \bf b}$, where for $i=1,\ldots,n$, $f_{x_i}$ denotes the partial derivative $\frac{\partial f}{\partial x_i}$. 
Morevover on ${\rm \bf b}$ the partial derivative $\frac{\partial f}{\partial v}$ does not vanish so we can decompose $\mathcal{B}$ into the disjoint union $\mathcal{B}^+ \sqcup \mathcal{B}^-$, where $\mathcal{B}^+$ (resp. $\mathcal{B}^-$) is the set of half-branches on which $\frac{\partial f}{\partial v} >0$ (resp. $\frac{\partial f}{\partial v} <0$).
This enables to define the following indices (Definition \ref{definitionlambda}):
$$\lambda^+= \sum_{{\rm \bf b} \in \mathcal{B}^+} \sigma({\rm \bf b}) \hbox{ and }\lambda^-= \sum_{{\rm \bf b} \in \mathcal{B}^-} \sigma({\rm \bf b}).$$
Then we define the following four indices (Definition \ref{definitiongamma}):
$$\begin{array}{l}
\gamma^{+,+}= \chi \big( f^{-1}(0) \cap \{x_1=a\} \cap B_\epsilon \big) - \chi \big( f^{-1}(\alpha) \cap \{x_1=a\} \cap B_\epsilon \big) , \cr
\gamma^{+,-}= \chi \big( f^{-1}(0) \cap \{x_1=-a\} \cap B_\epsilon \big) - \chi \big( f^{-1}(\alpha) \cap \{x_1=-a\} \cap B_\epsilon \big) ,\cr
\gamma^{-,+}= \chi \big( f^{-1}(0) \cap \{x_1=a\} \cap B_\epsilon \big) - \chi \big( f^{-1}(-\alpha) \cap \{x_1=a\} \cap B_\epsilon \big) ,\cr
\gamma^{-,-}= \chi \big( f^{-1}(0) \cap \{x_1=-a\} \cap B_\epsilon \big) - \chi \big( f^{-1}(-\alpha) \cap \{x_1=-a\} \cap B_\epsilon \big),
\end{array} $$
where $0< \alpha \ll a \ll \epsilon$. Our second generalization of the Khimshiashvili formula relates the Euler characteristic of the real Milnor fibres to these new indices. Namely in Theorem \ref{KhimPolar} we show that 
$$\chi \big( f^{-1}(-\delta) \cap B_\epsilon \big) =1-\lambda^- -\gamma^{-,-} = 1-\lambda^+ -\gamma^{-,+},$$
and that
$$\chi \big( f^{-1}(\delta) \cap B_\epsilon \big) =1-(-1)^n\lambda^- -\gamma^{+,+} = 1-(-1)^n \lambda^+ -\gamma^{+,-},$$
where $0 < \delta \ll \epsilon \ll 1$. Then we apply this result to the case where $\Sigma_f$ has dimension one.  In this case, we denote by $\mathcal{C}$ the set of  connected components of $\Sigma_f \setminus \{0\}$. For $v \in S^{n-1}$ generic, the function $v^*$ does not vanish on any half-branch ${\rm \bf c} \in \mathcal{C}$,  so we can decompose  $\mathcal{C}$ into the disjoint union $\mathcal{C}^+ \sqcup \mathcal{C}^-$, where $\mathcal{C}^+$ (resp. $\mathcal{C}^-$) is the set of half-branches on which $v^* >0$ (resp. $ v^*<0$). For each ${\rm \bf c} \in \mathcal{C}$, let  $\tau({\rm \bf c})$ be the value that the function $a \mapsto {\rm deg}_q \nabla f_{\vert x_1^{-1}(a)}$, $\{q\} = {\rm \bf c} \cap \{x_1 = a\}$, takes close to the origin.  Then we set $\gamma^+ = \sum_{{\rm \bf c} \in \mathcal{C}^+} \tau ({\rm \bf c})$ and $ \gamma^- = \sum_{{\rm \bf c} \in \mathcal{C}^-} \tau ({\rm \bf c})$. In this situation, Theorem \ref{KhimPolar} takes the following form (Theorem \ref{KhimOneDim}):
$$\chi \big( f^{-1}(-\delta) \cap B_\epsilon \big) = 1-\lambda^--\gamma^-=1-\lambda^++\gamma^+,$$
$$\chi \big( f^{-1}(\delta) \cap B_\epsilon \big) = 1-(-1)^n (\lambda^+-\gamma^-)=1-(-1)^n (\lambda^--\gamma^+),$$
where $0 < \delta \ll \epsilon \ll 1$.  Hence the indices $\lambda^+$, $\lambda^-$, $\gamma^+$ and $\gamma^-$ appear to be real versions of the first two L\^e numbers defined by Massey in \cite{MasseyInventiones1}. We note that the paper \cite{vanStratenWarmt} contains also formulas for the Euler characteristic of the real Milnor fibres of a function-germ with a one-dimensional critical locus.

In the complex case, the L\^e-Iomdine formula (\cite{LeIomdine,Iomdine}, see also \cite{MasseyInventiones1,MasseyIomdinePrepolar,MasseySiersma,SiersmaCommentarii1990} for improved versions) relates the Euler characteristic of the Milnor fibre of an analytic function-germ with one-dimensional singular set to the Milnor fibre of an analytic function-germ with an isolated singularity, given as the sum of the initial function and a sufficiently big power of a generic linear form.
As a corollary of Theorem \ref{KhimOneDim}, we establish a real version of this formula (Theorem \ref{RealLeIomdine}), i.e., a relation between the Euler characteristic of the real Milnor fibres of $f$ and the real Milnor fibres of a function of the type $f+ {v^*} ^k$, for $v \in S^{n-1}$ generic and $k \in \mathbb{N}$  big enough. 

We have also included some results on the local topology of locally closed definable sets. More precisely, we consider a locally closed definable set $X$ equipped with a Whitney stratification such that $0 \in X$, and a definable function $g : (X,0) \to (\mathbb{R},0)$ with an isolated critical point at the origin. In Lemma \ref{ArnoldWallStrat} we extend to this setting the results of Arnold and Wall mentioned above, i.e., we give  relations between the Euler characteristics of the sets $X \cap \{g \ ? \ \pm \delta \} \cap B_\epsilon$, where $0< \delta \ll \epsilon \ll 1$ and $? \in \{\le,\ge \}$, and the Euler characteristics of the sets $X \cap \{g \ ? \ 0 \} \cap S_\epsilon$, where $0<  \epsilon \ll 1$ and $? \in \{\le,\ge \}$. We give two corollaries (Corollaries \ref{Cor1ArnoldWallStrat} and \ref{LinkGenericForms}) when the stratum that contains $0$ has dimension greater than or equal to $1$.

The paper is organized as follows. In Section 2, we prove the first generalization of the Khimshiashvili formula based on Szafraniec's methods. In Section 3, we give the results on the local topology of locally closed definable sets. Section 4 contains the second generalization of the Khimshiashvili formula, based on the study of generic relative polar curves. In Section 5, we establish the real version of the L\^e-Iomdine formula. 

\vspace*{0,35cm}

{\it Acknowledgments.}
A large part of this paper was written during two visits of the author in the {\em Instituto de Ci\^{e}ncias Matem\'{a}ticas e de Computa\c{c}\~{a}o, Universidade de S\~{a}o Paulo - Campus de S\~{a}o Carlos}. The author thanks this institution, especially Maria Ruas and Nivaldo Grulha, for the financial support and the hospitality.  He also thanks Dirk Siersma  for fruitful discussions on one-dimensional singularities.

\section{Some general results on the real Milnor fibre}
Let $f:(\mathbb{R}^n, 0) \to (\mathbb{R},0)$ be a definable function-germ of class $C^r$, $r \ge 2$. By Lemma 10 in \cite{Bekka} or by the main theorem of \cite{TaLeLoi}, we can equip $f^{-1}(0)$ with a finite Whitney stratification that satisfies the Thom $(a_f)$-condition.

\begin{lemma}
There exists $\epsilon_0 >0$ such that for $0< \epsilon \le \epsilon_0$, there exists $\delta_\epsilon$ such that for $0 < \delta \le \delta_\epsilon$, the topological type of $f^{-1}(\delta) \cap B_\epsilon$ does not depend on the choice of the couple $(\epsilon,\delta)$.
\end{lemma}
\proof Let $\epsilon_0 >0$ be such that for $0 < \epsilon \le \epsilon_0$, the sphere $S_\epsilon$ intersects $f^{-1}(0)$ transversally. Then there exists a neighborhood $U_\epsilon$ of $0$ in $\mathbb{R}$ such that for each $\delta \in U_\epsilon$, the fibre $f^{-1}(\delta)$ intersects the sphere $S_\epsilon$ transversally. If it is not the case, then we can find a sequence of points $(p_m)_{m \in \mathbb{N}}$ in $S_\epsilon$ such that the vectors $\frac{p_m}{\vert p_m \vert}$ and $\frac{\nabla f (p_m)}{\vert \nabla f (p_m) \vert}$ are collinear, and such that the sequence converges to a point $p$ in $S_\epsilon \cap f^{-1}(0)$. If $S$ denotes the stratum of $f^{-1}(0)$ that contains $p$ then, applying the Thom $(a_f)$-condition, there exists a unit vector $V$ normal to $T_p S$ such that $\frac{p}{\vert p \vert}$ and $V$ are collinear. This contradicts the fact that $S_\epsilon$ intersects $f^{-1}(0)$ transversally. 

Now let us fix $\epsilon >0$ with $\epsilon \le \epsilon_0$. Let us choose $\delta_\epsilon >0$ such that $[0,\delta_\epsilon]$ is included in $U_\epsilon$ and $\delta$ is a regular value of $f$ for $0 < \delta \le \delta_\epsilon$. Let $(\epsilon_1,\delta_1)$ and $(\epsilon_2,\delta_2)$ be two couples with $0 < \epsilon_i \le \epsilon$ and $0< \delta_i \le \delta_{\epsilon_i}$ for $i=1,2$. If $\epsilon_1=\epsilon_2$ then the Thom-Mather first isotopy lemma implies that the fibres $f^{-1}(\delta_1) \cap B_{\epsilon_1}$ and $f^{-1}(\delta_2) \cap B_{\epsilon_2}$ are homeomorphic. Now assume that $\epsilon_1 < \epsilon_2$. By the same arguments as above, there exists a neighborhood $U$ of $0$ in $\mathbb{R}$ such that for each $\delta \not= 0$ in $U$, the distance function to the origin has no critical point on $f^{-1}(\delta) \cap (B_{\epsilon_2} \setminus \mathring{B_{\epsilon_1}})$. 
Let us choose $\delta_3 \not= 0$ in $U$ such that $0 < \delta_3 \le {\rm min}\{ \delta_1,\delta_2 \}$. By the first case, $f^{-1}(\delta_3) \cap B_{\epsilon_1}$ is homeomorphic to $f^{-1}(\delta_1) \cap B_{\epsilon_1}$ and
$f^{-1}(\delta_3) \cap B_{\epsilon_2}$ is homeomorphic to $f^{-1}(\delta_2) \cap B_{\epsilon_2}$. But, since the distance function to the origin has no critical points on $f^{-1}(\delta_3) \cap (B_{\epsilon_2} \setminus \mathring{B_{\epsilon_1}})$, the fibres $f^{-1}(\delta_3) \cap B_{\epsilon_2}$ and $f^{-1}(\delta_3) \cap B_{\epsilon_1}$ are homeomorphic.
\endproof
Of course a similar result is true for negative values of $f$. 
\begin{definition}
{\rm The (real) Milnor fibres of $f$ are the sets $f^{-1}(\delta) \cap B_\epsilon$ and $f^{-1}(-\delta) \cap B_\epsilon$, where $0 < \delta \ll \epsilon \ll 1$. }
\end{definition}
Sometimes we call $f^{-1}(\delta) \cap B_\epsilon$ (resp. $f^{-1}(-\delta) \cap B_\epsilon$) the positive (resp. negative) Milnor fibre of $f$. The Khimshiashvili formula \cite{Khimshiashvili} relates the Euler characteristic of the Milnor fibres to the topological degree of $\nabla f$ at the origin, when $f$ has an isolated singularity.

\begin{theorem}[The Khimshiashvili formula]\label{KhimshiashviliFormula} If $f$ has an isolated critical point at the origin then 
$$\chi \big( f^{-1}(\delta) \cap B_\epsilon \big) = 1-{\rm sign} (-\delta)^n {\rm deg}_0 \nabla f,$$
where $0 < \vert \delta \vert \ll \epsilon \ll 1$.
\end{theorem}
\proof We give a proof for we will need a similar argument later.   Let $U$ be a small open subset of $\mathbb{R}^n$ such that $0 \in U$ and $f$ is defined in $U$. We pertub $f$ in a Morse function $\tilde{f} : U \rightarrow \mathbb{R}$. Let $p_1,\ldots,p_k$ be the critical points of $\tilde{f}$, with respective indices $\lambda_1,\ldots,\lambda_k$. Let $\delta > 0$, by Morse theory we have:
$$\chi \big(f^{-1}([-\delta,\delta]) \cap B_\varepsilon\big) -\chi \big(f^{-1}(-\delta) \cap B_\varepsilon \big)=\sum_{i=1}^k (-1)^{\lambda_i}.$$
Actually we can choose $\tilde{f}$ sufficiently close to $f$ so that the $p_i$'s lie in $f^{-1}([-\frac{\delta}{4},\frac{\delta}{4}])$. Now the inclusion $f^{-1}(0) \cap B_\epsilon \subset f^{-1}([-\delta,\delta]) \cap B_\epsilon$ is a homotopy equivalence (Durfee \cite{Durfee} proved this result in the semi-algebraic case, but his argument holds in the $C^r$ definable case, see also \cite{CosteReguiat,Kurdyka98}) and $f^{-1}(0) \cap B_\epsilon$ is the cone over $f^{-1}(0) \cap S_\epsilon$, so $\chi \big( f^{-1}([-\delta,\delta]) \cap B_\epsilon \big)=1$. This gives the result for the negative Milnor fibre. To get the result for the positive one, it is enough to replace $f$ with $-f$. \endproof
The following formulas are due to Arnol'd \cite{Arnold} and Wall \cite{Wall}.
\begin{corollary}\label{ArnoldWall}
With the same hypothesis on $f$, we have:
$$\chi(\{f \le 0 \} \cap S_\varepsilon)= 1- \hbox{\em deg}_0 \nabla f, $$  $$\chi(\{f \ge 0 \} \cap S_\varepsilon)= 1+(-1)^{n-1}\hbox{\em deg}_0 \nabla f.$$
If $n$ is even,  we have:
$$\chi(\{f = 0 \} \cap S_\varepsilon )= 2-2\ \hbox{\em deg}_0 \nabla f.$$
\end{corollary}
\proof By a deformation argument due to Milnor \cite{Milnor}, $f(-\delta) \cap B_\varepsilon$, $\delta >0$, is homeomorphic to $\{ f \le -\delta \} \cap S_\varepsilon $, which is homeomorphic to $\{ f \le 0\} \cap S_\varepsilon $ if $\delta$ is very small. \endproof 

We start our study of the general case with an easy lemma.
\begin{lemma}\label{MilnorfibreLink}
Let $f: (\mathbb{R}^n,0) \to (\mathbb{R},0)$ be a definable function germ of class $C^r$, $r \ge 2$, and let $0 < \delta \ll \epsilon$. 
If $n$ is even then 
$$\chi \left( f^{-1}(-\delta) \cap B_\epsilon \right) = \chi \left( \{f \ge 0 \} \cap S_\epsilon \right),$$ and $$\chi \left( f^{-1}(\delta) \cap B_\epsilon \right) = \chi \left( \{f \le 0 \} \cap S_\epsilon \right).$$
If $n$ is odd then
$$\chi \left( f^{-1}(-\delta) \cap B_\epsilon \right) =2- \chi \left( \{f \ge 0 \} \cap S_\epsilon \right),$$ and $$\chi \left( f^{-1}(\delta) \cap B_\epsilon \right) =2- \chi \left( \{f \le 0 \} \cap S_\epsilon \right).$$
\end{lemma}
\proof If $n$ is even then $f^{-1}(-\delta) \cap B_\epsilon$ is an odd-dimensional manifold with boundary and so 
$$\chi \left( f^{-1}(-\delta) \cap B_\epsilon \right) =\frac{1}{2} \chi  \left( f^{-1}(-\delta) \cap S_\epsilon \right) =  \chi \left( \{ f \ge -\delta \} \cap S_\epsilon \right).$$
But for $\delta$ small, the inclusion $\{ f \ge 0 \} \cap S_\epsilon \subset\{ f \ge -\delta \} \cap S_\epsilon$  is a homotopy equivalence (see \cite{Durfee}).

If $n$ is odd then  $\{ f \ge -\delta \} \cap B_\epsilon$ is an odd-dimensional manifold with corners. Rounding the corners, we get 
$$\displaylines{ \quad \chi \left( \{ f \ge -\delta \} \cap B_\epsilon \right) =\frac{1}{2} \Big(  \chi \left( f^{-1}(-\delta) \cap B_\epsilon \right) +  \chi  \left( \{ f \ge -\delta \} \cap S_\epsilon \right)  \hfill  \cr
\hfill - \chi  \left( f^{-1}(-\delta) \cap S_\epsilon \right) \Big)
= \frac{1}{2} \left(  \chi \left( f^{-1}(-\delta) \cap B_\epsilon \right) +  \chi  \left(\{ f \ge -\delta \} \cap S_\epsilon \right) \right). \quad \cr
}$$
 But the inclusion $ \{ f \ge 0\} \cap B_\epsilon  \subset  \{ f \ge -\delta \} \cap B_\epsilon  $ is a homotopy equivalence and so 
$$\chi \left( \{ f \ge -\delta \} \cap B_\epsilon \right) =1.$$ \endproof
For the rest of this section, we assume that the structure is polynomially bounded. The technics developed and the results proved by Szafraniec \cite{SzafraniecTopology86} (see also \cite{Bruce89}) are valid in this context. Let $\omega(x)=x_1^2+\cdots+x_n^2$. Then there exists an integer $d>0$ sufficiently big such that $g_+=f - \omega^d$ and $g_-=-f-\omega^d$ have an isolated critical point at the origin. 
Moreover Szafraniec showed that 
$$\chi \left( \{f \le 0 \} \cap S_\epsilon \right) =1-{\rm deg}_0 \nabla g_+ \hbox{ and }  \chi \left( \{f \ge 0 \} \cap S_\epsilon \right) =1- {\rm deg}_0 \nabla g_-.$$
Applying the previous lemma, we can state our first generalization of the Khimshiashvili formula. 
\begin{corollary}\label{KhimPolyBounded}
If $0< \delta \ll \epsilon$, we have:
$$\chi \left( f^{-1}(-\delta) \cap B_\epsilon \right) =  1-(-1)^n {\rm deg}_0 \nabla g_-,$$
and 
$$\chi \left( f^{-1}(\delta) \cap B_\epsilon \right) =  1-(-1)^n {\rm deg}_0 \nabla g_+.$$
\end{corollary}

In general, the exponent $d$ is difficult to estimate. However, in case of a weighted-homogeneous polynomial, Szafraniec \cite{SzafraniecGlasgow91} provided another method which is completely effective.

Let $f:\mathbb{R}^{n} \to \mathbb{R}$ be a real weighted homogeneous polynomial function of type $(d_{1},\cdots,d_{n};d)$ with $\nabla f(0)=0.$ 
Let $p$ be the smallest positive integer such that $2p>d$ and each $d_{i}$ divides $p.$ Also denote by $a_{i}=\displaystyle{\frac{p}{d_{i}}}$ and

$$\omega=\frac{x_{1}^{2a_{1}}}{2a_{1}}+\cdots+\frac{x_{n}^{2a_{n}}}{2a_{n}}.$$
Now consider $g_{1}=f-\omega$ and $g_{2}=-f-\omega.$  Szafraniec proved that $g_1$ and $g_2$ have an algebraically isolated critical point at the origin and that 
$$\chi \left( \{f \le 0 \} \cap S^{n-1} \right) =1-{\rm deg}_0 \nabla g_1 \hbox{ and }  \chi \left( \{f \ge 0 \} \cap S^{n-1} \right) =1- {\rm deg}_0 \nabla g_2.$$
Applying Lemma \ref{MilnorfibreLink}, we obtain the following Khimshiashvili's type formula for the fibres of a real weighted homogeneous polynomial.
\begin{corollary}\label{FibreMilnorHomogeneous}
We have 
$$\chi \left( f^{-1}(-1) \right) =  1-(-1)^n {\rm deg}_0 \nabla g_2,$$
and 
$$\chi \left( f^{-1}(1)  \right) =  1-(-1)^n {\rm deg}_0 \nabla g_1.$$
\end{corollary}
Note that ${\rm deg}_0 \nabla g_1$ and ${\rm deg}_0 \nabla g_2$ can be computed algebraically thanks to the Eisenbud-Levine-Khimshiashvili formula \cite{EL,Khimshiashvili} because they have an algebraically isolated zero at the origin.

Let us apply this corollary to the examples presented in \cite{SzafraniecGlasgow91}.
\begin{enumerate}
\item Let $f(x,y,z)=x^2y-y^4-yz^3$. By \cite{SzafraniecGlasgow91}, we have that 
${\rm deg}_0 \nabla g_1={\rm deg}_0 \nabla g_2=1$. So $\chi \left( f^{-1}(-1) \right)=\chi \left( f^{-1}(1) \right)=2$.
\item Let $f(x,y,z)=x^3+x^2z-y^2$. By \cite{SzafraniecGlasgow91}, we have that 
${\rm deg}_0 \nabla g_1=1$ and ${\rm deg}_0 \nabla g_2=-1$. So $\chi \left( f^{-1}(-1) \right)=0$ and $\chi \left( f^{-1}(1) \right)=2$.
\item Let $f(x,y,z)=x^3-xy^2+xyz+2x^2y -2y^3-y^2z-xz^2+yz^2$. Then by \cite{SzafraniecGlasgow91}, 
${\rm deg}_0 \nabla g_1=3$ , so $\chi \left( f^{-1}(1) \right)=4$.
\end{enumerate}

\section{Some results on the topology of locally closed definable sets}
Let $X$ be a locally closed definable set. We assume that $0$ belongs to $X$.
We equip $X$ with a finite definable $C^r$, $r\ge 2$, Whitney stratification. The fact that such a stratification exists is due to Loi \cite{Loi98}. Recently Nguyen, Trivedi and Trotman \cite{NguyenTrivediTrotman} gave another proof of this result. We denote by $S_0$ the stratum that contains $0$. 

Let $g : (X,0) \to (\mathbb{R},0)$ be a definable function that is the restriction to $X$ of a definable function $G$ of class $C^r$, $r \ge 2$, defined in a neighborhood of the origin. We assume that $g$ has at worst an isolated critical point (in the stratified sense) at the origin. As in the previous section, the positive and the negative real Milnor fibres of $g$ are the sets $g^{-1}(\delta) \cap X \cap B_\epsilon$ and $g^{-1}(-\delta) \cap X \cap B_\epsilon$, where $0< \delta \ll \epsilon \ll 1$.

\begin{lemma}\label{ArnoldWallStrat} For $0 < \delta \ll \epsilon \ll 1$, we have 
$$\chi \big( X \cap g^{-1}(-\delta) \cap B_\epsilon \Big)= 
\chi \big( X \cap \{ g \le 0 \} \cap S_\epsilon \Big),$$
and 
$$\chi \big( X \cap g^{-1}(\delta) \cap B_\epsilon \Big)= 
\chi \big( X \cap \{ g \ge 0 \} \cap S_\epsilon \Big).$$
\end{lemma}
\proof Using the methods developed in \cite{DutertreJofSingProcTrot}, we can assume that the critical points of $g$ on $X \cap S_\epsilon$ are isolated, that they lie in $\{ g \not= 0 \}$ and that they are outwards-pointing (resp. inwards-pointing) in $\{ g > 0 \}$ (resp. $\{g < 0 \}$). Let us denote them by $\{p_1,\ldots,p_s\}$. 

We recall that if $Z \subset \mathbb{R}^n$ is a locally closed definable set, equipped with a Whitney stratification and $p$ is an isolated critical point of a definable function $\phi : Z \rightarrow \mathbb{R}$, restriction to $Z$ of a $C^2$-definable function $\Phi$, then the index of $\phi$ at $p$ is defined as follows:
$${\rm ind}(\phi,Z,p)= 1 - \chi \big( Z \cap \{ \phi = \phi(p)-\eta \} \cap B_r(p) \big),$$
where $0< \eta \ll r \ll 1$ and $B_r(p)$ is the closed ball of radius $r$ centered at $p$. 

As in \cite{DutertreJofSingProcTrot}, Section 3, we can apply the results proved in \cite{DutertreManuscripta12}. Namely, by Theorem 3.1 in \cite{DutertreManuscripta12}, we can write 
$$\chi \big( \{ g \le 0 \} \cap X \cap S_\epsilon \big) = \sum_{i \ \vert \ g(p_i) < 0 } {\rm ind}(g,X \cap S_\epsilon,p_i),$$
and for $0< \delta \ll \epsilon$,
$$\chi \big( \{ g \le \delta \} \cap X \cap B_\epsilon \big) = \sum_{i \ \vert \ g(p_i) < 0 } {\rm ind}(g,X \cap B_\epsilon,p_i) + {\rm ind}(g,X,0).$$
By Lemma 2.1 in \cite{DutertreManuscripta12}, ${\rm ind}(g, X \cap S_\epsilon, p_i)={\rm ind}(g,X \cap B_\epsilon,p_i)$ if $g(p_i)<0$. Moreover, ${\rm ind}(g,X,0) =1- \chi \big( g^{-1}(-\delta) \cap X \cap B_\epsilon \big)$ and, as explained in the proof of Theorem \ref{KhimshiashviliFormula}, $\chi \left( \{g \le \delta \} \cap X \cap B_\epsilon \right)=1$ if $\delta$ is small enough. 
Combining these observations, we find that 
$$\chi \big( X \cap g^{-1}(-\delta) \cap B_\epsilon \big) =
\chi \big( X \cap \{ g \le 0 \} \cap S_\epsilon \big).$$ \endproof

\begin{remark} {\rm We believe that it is possible to establish these equalities applying a stratified version of the Milnor deformation argument mentionned in the proof of Corollary \ref{ArnoldWall}. This is done by Comte and Merle in \cite{ComteMerle} when $X$ is conic and $g$ is a generic linear form.}
\end{remark}

For the rest of this section, we will denote by ${\rm Lk}(Y)$ the link at the origin of a definable set $Y$.

\begin{corollary}\label{Cor1ArnoldWallStrat} Assume that ${\rm dim}\ S_0 >0$ and that $g_{\vert S_0}$ has no critical point at $0$, i.e., $g^{-1}(0)$ intersects $S_0$ transversally at $0$. Then the following equalities hold:
$$\chi \big( {\rm Lk}(X \cap \{g \le 0 \}) \big)= \chi \big( {\rm Lk}(X \cap \{g \ge 0 \}) \big) =1,$$
and 
$$\chi \big( {\rm Lk}(X  ) \big)+ \chi \big( {\rm Lk}(X \cap \{g = 0 \}) \big)=2.$$
\end{corollary}
\proof If $g_{\vert S_0}$ has no critical point at $0$, then $g : X \to \mathbb{R}$ is a stratified submersion in a neighborhood of $0$. Furthermore for $0 < \epsilon \ll 1$, the sphere $S_\epsilon$ intersects $X \cap \{g=0\}$ transversally, so $0$ is a regular value of $g_{\vert X \cap B_\epsilon}$. Therefore if $\delta$ is small enough, 
$$\chi \big( X \cap \{g=-\delta\} \cap B_\epsilon \big)=
\chi \big( X \cap \{g=\delta\} \cap B_\epsilon \big)=
\chi \big( X \cap \{g=0\} \cap B_\epsilon \big)=1.$$
It is enough to apply the previous lemma and then the Mayer-Vietoris sequence. \endproof

For $v \in S^{n-1}$, we denote by $v^*$ the function $v^*(x)= \langle v, x \rangle$, where $\langle \ , \ \rangle$ is the standard scalar product. The previous corollary applies to a generic linear form $v^*$.
\begin{corollary}\label{LinkGenericForms}
Assume that ${\rm dim}\ S_0 >0$. If $v \notin S^{n-1} \cap (T_0 S_0)^\perp$, then 
$$\chi \big( {\rm Lk}(X \cap \{ v^* \le 0 \} ) \big) = 
\chi \big( {\rm Lk}(X \cap \{ v^* \ge 0 \} ) \big)=1,$$
and 
$$\chi \big( {\rm Lk}(X  ) \big)+ \chi \big( {\rm Lk}(X \cap \{v^* = 0 \}) \big)=2.$$
\end{corollary}
\proof If $v \notin (T_0 S_0)^\perp$, then $v^*_{\vert S_0}$ has no critical point at $0$. \endproof
Let us relate this corollary to results that we proved in earlier papers. Combining Theorem 5.1 in \cite{DutertreJofSingProcTrot} and the comments after Theorem 2.6 in \cite{DutertreIsrael}, we can write that if ${\rm dim}\ S_0 >0$, 
$$\chi \big( {\rm Lk}(X  ) \big) +\frac{1}{g_n^{n-1}} \int_{G_n^{n-1}} \chi \big( {\rm Lk}(X \cap H) \big) dH =2,$$
where $G_n^{n-1}$ is the Grassmann manifold of linear hyperplanes in $\mathbb{R}^n$ and $g_n^{n-1}$ is its volume. This last equality is based on the study of the local behaviour of the generalized Lipschitz-Killing curvatures made in \cite{ComteMerle} and \cite{DutertreJofSingProcTrot}. We see that it is actually a direct consequence of Corollary \ref{LinkGenericForms}, which gives a more precise result on the local topology of locally closed definable sets. Similarly for $0< k < {\rm dim}\ S_0$, we know that 
$$-\frac{1}{g_n^{n-k-1}} \int_{G_n^{n-k-1}} \chi \big( {\rm Lk}(X \cap H) \big) dH + \frac{1}{g_n^{n-k+1}} \int_{G_n^{n-k+1}} \chi \big( {\rm Lk}(X \cap L) \big) dL = 0,$$
where $G_n^{n-k}$ is the Grassman manifold of $k$-dimensional vector spaces in $\mathbb{R}^n$ and $g_n^{n-k}$ is its volume. In fact a recursive application of Corollary \ref{LinkGenericForms} shows that $\chi \big( {\rm Lk}(X \cap H) \Big)=\chi \big( {\rm Lk}(X \cap L) \Big)$ for $H$ generic in $G_n^{n-k-1}$ and $L$ generic in $G_n^{n-k+1}$. 

Let us give another application of Corollary \ref{LinkGenericForms} to the topology of real Milnor fibres. As in the previous section, $f : (\mathbb{R}^n,0) \to (\mathbb{R},0)$ is the germ at the origin of a definable function of class $C^r$, $r \ge 2$. We assume that $f^{-1}(0)$ is equipped with a  finite Whitney stratification that satisfies the Thom $(a_f)$-condition. Let $S_0$ be the stratum that contains $0$.
\begin{corollary}\label{Corollary3.5}
If ${\rm dim}\ S_0 >0$ and if $v \notin S^{n-1} \cap (T_0 S_0)^\perp$, then for $0 < \delta \ll \epsilon \ll 1$, we have 
$$\chi \big( f^{-1}(\delta) \cap B_\epsilon \big)=
\chi \big( f^{-1}(\delta) \cap \{ v^*=0 \} \cap B_\epsilon \big),$$
and
$$\chi \big( f^{-1}(-\delta) \cap B_\epsilon \big)=
\chi \big( f^{-1}(-\delta) \cap \{ v^*=0 \} \cap B_\epsilon \big).$$
\end{corollary}
\proof Applying Corollary \ref{LinkGenericForms} to the sets $\{f \ge 0 \}$ and $\{f \le 0 \}$, we get that
$$\chi \big( {\rm Lk}( \{ f  \ ? \ 0 \}  ) \big)+ \chi \big( {\rm Lk}(\{ f \ ? \ 0 \}  \cap \{v^* = 0 \}) \big)=2,$$
where $? \in \{ \le, \ge \}$. Lemma \ref{MilnorfibreLink} applied to $f$ and $f_{\vert \{ v^*=0 \}}$ gives the result. \endproof
In the next section, we will give a generalization of this result based on generic relative polar curves.

\section{Milnor fibres and relative polar curves}
Let $f :(\mathbb{R}^n,0) \to (\mathbb{R},0)$ be a definable function-germ of class $C^r$, $r \ge 2$. We will give a second generalization of the Khimshiashvili formula in this setting. For this we need first to study the behaviour of a generic linear function on the fibres of $f$ and the behaviour of $f$ on the fibres of a generic linear function. 

We start with a study of the critical points of $v^*_{\vert f^{-1}(\delta)}$ for $\delta$ small and $v$ generic in $S^{n-1}$. Let 
$$\Gamma_v =\left\{ x \in \mathbb{R}^n \setminus \Sigma_f \ \vert \ {\rm rank}(\nabla f(x), v ) <2 \right\}.$$

We will need a first genericity condition. We can equip $f^{-1}(0)$ with a finite Whitney stratification that satisfies the Thom $(a_f)$-condition.
\begin{lemma}\label{Lemma4.1}
There exists a definable set $\Sigma_1 \subset S^{n-1}$ of positive codimension such that if $v \notin \Sigma_1$, then $\{ v^* = 0\}$ intersects $f^{-1}(0) \setminus \{0\}$ transversally in a neighborhood of the origin. 
\end{lemma}
\proof It is a particular case of Lemma 3.8 in \cite{DutertreGeoDedicata2012}. \endproof

\begin{lemma}\label{Lemma4.2}
If $v \notin \Sigma_1$ then $\Gamma_v \cap f^{-1}(0) = \emptyset$.
\end{lemma}
\proof If it is not the case then we can find an arc $\alpha : [0, \nu[ \to f^{-1}(0)$ such that $\alpha (0)=0$ and for $0< s < \nu$, $\nabla f(\alpha(s)) \not= 0$ and ${\rm rank}(\nabla f (\alpha(s)), v) < 2$. Let $S$ be the stratum that contains $\alpha (]0,\nu [)$. Since $\nabla f (\alpha(s))$ is normal to $T_{\alpha(s)} S$, the points in $\alpha(]0,\nu[)$ are critical points of $v^*_{\vert S}$ and hence lie in $\{ v^*=0\}$. This contradicts Lemma \ref{Lemma4.1}. \endproof

\begin{corollary}\label{Corollary4.3}
If $v \notin \Sigma_1$ then $\Gamma_v \cap \{v^* =  0 \} = \emptyset$.
\end{corollary}
\proof As in the proof of the previous lemma, we see that if $\Gamma_v \cap \{v^*=0\} \not= \emptyset$ then $\Gamma_v \cap f^{-1}(0) \not= \emptyset$. \endproof

\begin{lemma}\label{Lemma4.4}
There exists a definable set $\Sigma_2 \subset S^{n-1}$ of positive codimension such that if $v \notin \Sigma_2$, $\Gamma_v$ is a curve (possibly empty) in the neighbourhood of the origin.
\end{lemma}
\proof Let 
$$M = \left\{ (x,y) \in \mathbb{R}^n \times \mathbb{R}^n \ \vert \ {\rm rank}(\nabla f(x),y) < 2 \right\}.$$
Let $p=(x_0,y_0)$ be a point in $M \setminus (\Sigma_f \times \mathbb{R}^n)$. We can assume that $f_{x_1}(x_0) \not= 0$. Therefore locally $M \setminus (\Sigma_f \times \mathbb{R}^n)$ is given by the equations $m_{12}(x,y)=\cdots=m_{1n}(x,y)=0$, where 
$$m_{ij}(x,y)= \left\vert \begin{array}{cc}
f_{x_i}(x) & f_{x_j}(x) \cr
y_i & y_j \cr
\end{array} \right\vert.$$
The Jacobian matrix of the mapping $(m_{12},\ldots,m_{1n})$ has the following form
$$\left( \begin{array}{ccccccc}
* & \cdots & * & -f _{x_2} & f_{x_1} &  \cdots & 0 \cr
\vdots & \ddots & \vdots & \vdots &  \vdots & \ddots & \vdots \cr
* & \cdots & * & -f _{x_n} &  0 & \cdots & f_{x_1} \cr
\end{array}
\right).$$
This implies that $M \setminus (\Sigma_f \times \mathbb{R}^n)$ is a $C^{r-1}$ manifold of dimension $n+1$. The Bertini-Sard theorem (\cite{BochnakCosteRoy}, 9.5.2) implies that the discriminant $D$ of the projection
$$\begin{array}{ccccc}
\pi_y & : & M \setminus (\Sigma_f \times \mathbb{R}^n) & \to & \mathbb{R}^n \cr
 & & (x,y) & \mapsto & y \cr
\end{array}$$
is a definable set of dimension less than or equal to $n-1$. Hence for all $v \in S^{n-1} \setminus D$, the dimension of $\pi_y^{-1}(v)$ is less than or equal to $1$. But $\pi_y^{-1}(v)$ is exactly $\Gamma_v$ and we set $\Sigma_2=D \cap S^{n-1}$. \endproof

\begin{corollary}\label{Corollary4.5}
Let $v \in S^{n-1}$ be such that $v \notin \Sigma_2$. There exists $\delta_v'$ such that for $0 < \vert \delta \vert \le \delta_v'$, the critical points of $v^*_{\vert f^{-1}(\delta)}$ are Morse critical points in a neighborhood of the origin. 
\end{corollary}
\proof After a change of coordinates, we can assume that $v=e_1=(1,0,\ldots,0) \in \mathbb{R}^n$ and so that $v^*(x)=x_1$. 

Let $p$ be a point in $\Gamma_v=\Gamma_{e_1}$. If $f_{x_1}(p) = 0$ then, since the minors $m_{1i} = \frac{\partial(f,x_1)}{\partial (x_1,x_i)}$, $i=2,\ldots,n$, vanish at $p$, $f_{x_i}(p) =0$ for $i=2,\ldots,n$ and so $p$ belongs to $\Sigma_f$, which is impossible. Therefore $f_{x_1}(p) \not= 0$ and by the proof of Lemma \ref{Lemma4.4}, we conclude that $\Gamma_{e_1}$ is defined by the vanishing of the minors $m_{1i}$, $i=2,\ldots,n$, and that 
$${\rm rank} \big( \nabla m_{12},\ldots,m_{1n} \big) = n-1$$
along $\Gamma_{e_1}$. Let ${\rm \bf a}$ be an arc (i.e., a connected component) of $\Gamma_{e_1}$, and let $\alpha : [0,\nu[ \to \bar{{\rm \bf a}}$ be a $C^r$ definable parametrization such that $\alpha(0)=0$ and $\alpha(]0,\nu[) \subset {\rm \bf a}$. Since $f$ does not vanish on ${\rm \bf a}$, the function $f \circ \alpha$ is strictly monotone which implies that for $s \in ]0,\nu[$, $(f \circ \alpha)'(s) = \langle \nabla f (\alpha(s)),\alpha'(s) \rangle \not=0$. Hence the vectors $$\nabla f(\alpha(s)), \nabla m_{12}(\alpha(s)) ,\ldots, \nabla m_{1n}(\alpha(s))$$ are linearly independent since the $\nabla m_{1i}(\alpha(s))$'s are orthogonal to $\alpha'(s)$. By Lemma 3.2 in \cite{SzafraniecManus94}, this is equivalent to the fact that the function $x_1 : f^{-1}(f(\alpha(s)) \to \mathbb{R}$ has a non-degenerate critical point at $\alpha(s)$. It is easy to conclude because $\Gamma_{e_1}$ has a finite numbers of arcs. \endproof

From now on, we will work with $v \in S^{n-1}$ such that $v \notin \Sigma_1 \cup \Sigma_2$. After a change of coordinates, we can assume that $v =e_1=(1,0,\ldots,0)$ and so the conclusions of Lemma \ref{Lemma4.1}, Lemma \ref{Lemma4.2}, Corollary \ref{Corollary4.3}, Lemma \ref{Lemma4.4} and Corollary \ref{Corollary4.5} are valid for $\Gamma_{x_1}$ and $\{x_1 = 0 \}$. Let us study the points of $\Gamma_{x_1}$ more accurately. By the previous results, we know that if $p$ is a point of $\Gamma_{x_1}$ close to the origin then $p$ is a Morse critical point of ${x_1}_{\vert f^{-1}(f(p))}$, $f_{x_1}(p) \not= 0$, $x_1(p) \not=0$ and $f(p) \not=0$. 

\begin{lemma}\label{Lemma4.6} 
Let $p$ be a point in $\Gamma_{x_1}$ close to the origin. Let $\mu(p)$ be the Morse index of ${x_1}_{\vert f^{-1}(f(p))}$ at $p$. Then $p$ is a Morse critical point of $f_{\vert x_1^{-1}(x_1(p))}$ and if $\theta(p)$ is the Morse index of $f_{\vert x_1^{-1}(x_1(p))}$ at $p$ then 
$$(-1)^{\mu(p)} =(-1)^{n-1} {\rm sign} (f_{x_1}(p))^{n-1} (-1)^{\theta(p)}.$$
\end{lemma}
\proof By Lemma 3.2 in \cite{SzafraniecManus94}, we know that 
$${\rm det} \big[\nabla f(p),\nabla f_{x_2}(p),\ldots,\nabla f_{x_n}(p) \big] \not= 0,$$
and that 
$$(-1)^{\mu(p)} =(-1)^{n-1} {\rm sign} (f_{x_1}(p))^{n} {\rm sign }\ {\rm det} \big[\nabla f(p),\nabla f_{x_2}(p),\ldots,\nabla f_{x_n}(p) \big].$$
But $\nabla f(p)= f_{x_1}(p) e_1$ and so ${\rm det} \big[ e_1,\nabla f_{x_2}(p),\ldots,\nabla f_{x_n}(p) \big]\not =0$ and 
$$(-1)^{\mu(p)} = (-1)^{n-1} {\rm sign} (f_{x_1}(p))^{n-1} {\rm sign }\ {\rm det} \big[e_1,\nabla f_{x_2}(p),\ldots,\nabla f_{x_n}(p) \big].$$
Still using Lemma 3.2 in \cite{SzafraniecManus94}, we see that $p$ is a Morse critical point of $f_{\vert x_1^{-1}(x_1(p))}$ and that 
$$(-1)^{\mu(p)} =(-1)^{n-1} {\rm sign} (f_{x_1}(p))^{n-1} (-1)^{\theta(p)}.$$ \endproof

\begin{lemma}\label{Lemma4.7}
Let $p$ be a point in $\Gamma_{x_1}$ close to the origin. Then $${\rm det} \big[ \nabla f_{x_1}(p), \nabla f_{x_2}(p),\ldots, \nabla f_{x_n}(p) \big] \not=0$$ and 
$$(-1)^{\theta(p)} = {\rm sign} \left( x_1(p) f_{x_1}(p) \right) {\rm sign} \ {\rm det} \big[ \nabla f_{x_1}(p), \nabla f_{x_2}(p),\ldots, \nabla f_{x_n}(p) \big] .$$
\end{lemma}
\proof Since ${\rm det} \big[e_1, \nabla f_{x_2}(p),\ldots, \nabla f_{x_n}(p) \big] \not=0$, we can write 
$$\nabla f_{x_1}(p) = \beta(p) e_1 + \sum_{i=2}^n \beta_i(p) \nabla f_{x_i}(p),$$
and so, 
$${\rm det} \big[ \nabla f_{x_1}(p), \nabla f_{x_2}(p),\ldots, \nabla f_{x_n}(p) \big] = \beta (p) \big[ e_1, \nabla f_{x_2}(p),\ldots, \nabla f_{x_n}(p) \big] .$$
Let $\alpha : [0,\nu[ \to \overline{\Gamma_{x_1}}$ be a parametrization of the arc that contains $p$. We have 
$$ ( f_{x_1} \circ \alpha )'(s)=\langle \nabla f_{x_1}(\alpha(s)) , \alpha'(s) \rangle = \beta(\alpha(s)) \langle e_1, \alpha' (s) \rangle = \beta(\alpha(s)) (x_1 \circ \alpha)'(s).$$
But since $f_{x_1}$ and $x_1$ do not vanish on $\Gamma_{x_1}$, $(f_{x_1} \circ \alpha)'(s)$ and $(x_1 \circ \alpha)'(s)$ do not vanish for $s$ small. Therefore for $p$ close to the origin, $\beta (p) \not=0$ and 
$${\rm sign}\ \beta(p)= {\rm sign} \left( x_1(p) f_{x_1}(p) \right) .$$
\endproof

Let $\mathcal{B}$ be the set of connected components of $\Gamma_{x_1}$. If ${\rm \bf b} \in \mathcal{B}$ then ${\rm \bf b}$ is a half-branch on which the functions $f_{x_1}$ and ${\rm det}\big[ \nabla f_{x_1},\ldots,\nabla f_{x_n} \big]$ have constant sign. So we can decompose $\mathcal{B}$ into the disjoint union $\mathcal{B}^+ \sqcup \mathcal{B}^-$ where $\mathcal{B}^+$ (resp. $\mathcal{B}^-$) is the set of half-branches on which $f_{x_1} >0$ (resp. $f_{x_1}<0$). If ${\rm \bf b} \in \mathcal{B}$, we denote by $\sigma({\rm \bf b})$ the sign of ${\rm det}\big[ \nabla f_{x_1},\ldots,\nabla f_{x_n} \big]$ on ${\rm \bf b}$. 
\begin{definition}\label{definitionlambda}
{\rm We set $\lambda^+= \sum_{{\rm \bf b} \in \mathcal{B}^+} \sigma({\rm \bf b})$ and $\lambda^-= \sum_{{\rm \bf b} \in \mathcal{B}^-} \sigma({\rm \bf b})$.}
\end{definition}
\begin{remark}
{\rm If $f$ has an isolated critical point at the origin then for $\eta \not=0$ small enough, $(\nabla f)^{-1}(\eta,0,\ldots,0)$ is exactly $\Gamma_{x_1} \cap f_{x_1}^{-1}(\eta)$. Moreover if $p \in \Gamma_{x_1} \cap f_{x_1}^{-1}(\eta)$, then ${\rm sign}\ {\rm det} \big[ \nabla f_{x_1}(p), \nabla f_{x_2}(p),\ldots, \nabla f_{x_n}(p) \big]  \not=0$. Hence $(\eta,0,\ldots,0)$ is a regular value of $\nabla f$ and so
$${\rm deg}_0 \nabla f =\sum_{p \in \Gamma_{x_1} \cap f_{x_1}^{-1}(\eta) } {\rm sign}\ {\rm det} \big[ \nabla f_{x_1}(p), \nabla f_{x_2}(p),\ldots, \nabla f_{x_n}(p) \big].$$
If $\eta >0$ (resp. $\eta <0$), this implies that ${\rm deg}_0 \nabla f = \lambda^+$ (resp. $\lambda^-$).}
\end{remark}

The following lemma will enable us to define other indices associated with $f$ and $x_1$. 
\begin{lemma}\label{Lemma4.10} 
There exists $\epsilon_0 >0$ such that for $0< \epsilon \le \epsilon_0$, there exists $a_\epsilon >0$ such that for $0 < a \le a_\epsilon$, there exists $\alpha_{a,\epsilon}>0$ such that for $0 <\alpha \le \alpha_{a,\epsilon}$, the topological type of $f^{-1}(\alpha) \cap \{x_1 = a \} \cap B_\epsilon$ does not depend on the choice of the triplet $(\epsilon,a,\alpha)$.
\end{lemma}
\proof For $a>0$ small enough, we define $\beta(a)$ by 
$$\beta (a)= {\rm inf} \big\{ \vert f (p) \vert \ \vert \ p \in \Gamma_{x_1} \cap \{  x_1  = a \} \big\}.$$
The function $\beta$ is well defined because $\Gamma_{x_1} \cap \{x_1 = a \} $ is finite and $\beta (a) >0$. Moreover it is definable and so it is continuous on a small interval of the form $]0,u[$. This implies that the set 
$$\mathcal{O} = \left\{ (a,\alpha) \in \mathbb{R} \times \mathbb{R}^* \ \vert \ a \in ]0,u[ \hbox{ and }
0 < \alpha < \beta (a) \right\}$$ is open and connected.

Since $\{x_1=0\}$ intersects $f^{-1}(0) \setminus \{0\}$ transversally (in the stratified sense), $\{x_1=0\} \cap f^{-1}(0) \setminus \{0\}$ is Whitney stratified, the strata being the intersections of $\{x_1=0\}$ with the strata of $f^{-1}(0) \setminus \{0\}$. 

Let $\epsilon_0 >0$ be such that for $0 < \epsilon \le \epsilon_0$, the sphere $S_\epsilon$ intersects $\{x_1=0\} \cap f^{-1}(0)$ transversally. Then there exists a neighborhood $\mathcal{U}_\epsilon$ of $(0,0)$ in $\mathbb{R}^2$ such that for each $(a,\alpha)$ in $(\mathbb{R} \times \mathbb{R}^*) \cap \mathcal{U}_\epsilon$, the fibre $f^{-1}(\alpha) \cap \{x_1 = a \}$ intersects $S_\epsilon$ transversally. If it is not the case, then we can find a sequence of points $(p_m)_{m \in \mathbb{N}}$ in $S_\epsilon$ such that the vectors $e_1$, $\frac{p_m}{\vert p_m \vert}$ and $\frac{\nabla f (p_m) }{\vert \nabla (p_m) \vert}$ are linearly dependent, and such that the sequence converges to a point $p$ in $S_\epsilon \cap f^{-1}(0) \cap \{x_1 =0\}$. If $S$ denotes the stratum of $f^{-1}(0)$ that contains $p$ then, applying the Thom $(a_f)$-condition and the method of Lemma 3.7 in  
\cite{DutertreJofSingProcTrot}, there exists a unit vector $v$ normal to $T_p S$ such that the vectors $e_1$, $\frac{p}{\vert p \vert}$ and $v$ are linearly dependent. But $e_1$ and $v$ are linearly independent for $\{x_1=0\}$ intersects $S$ transversally at $p$. Therefore $S_\epsilon$ does not intersect $S \cap \{x_1 = 0\}$ transversally at $p$, which is a contradiction. Moreover we can assume that $\mathcal{U}_\epsilon \cap \mathcal{O}$ is connected. 

Now let us fix $\epsilon >0$ with $\epsilon \le \epsilon_0$. Let us choose $a_\epsilon >0$ such that $a_\epsilon <u$ and the interval $]0,a_\epsilon]$ is included in $\mathcal{U}_\epsilon$. For each $a \in ]0,a_\epsilon]$, there exists $\alpha'_{a,\epsilon}$ such that $ \{a \} \times ]0,\alpha'_{a,\epsilon}]$ is included in $\mathcal{U}_\epsilon$. We choose $\alpha_{a,\epsilon}$ such that $\alpha_{a,\epsilon} \le \alpha'_{a,\epsilon}$ and $\alpha_{a,\epsilon} < \beta (a)$, which implies that $(a,\alpha)$ is a regular value of $(x_1,f)$ for $0 < \alpha \le \alpha_{a,\epsilon}$. 

Let $(\epsilon_1,a_1,\alpha_1)$ and $(\epsilon_2,a_2,\alpha_2)$ be two triplets with $0 < \epsilon_i \le \epsilon$, $0<a_i \le a_{\epsilon_i}$ and $0 < \alpha_i \le \alpha_{a_i,\epsilon_i}$ for $i=1,2$. If $\epsilon_1=\epsilon_2$ then the Thom-Mather first isotopy lemma implies that the fibres $f^{-1}(\alpha_1) \cap \{x_1=a_1 \} \cap B_{\epsilon_1}$ and $f^{-1}(\alpha_2) \cap \{x_1=a_2 \} \cap B_{\epsilon_2}$ are homeomorphic, because $(a_1,\alpha_1)$ and $(a_2,\alpha_2)$ belong to the connected set $\mathcal{U}_{\epsilon_1} \cap \mathcal{O}$.

Now assume that $\epsilon_1 < \epsilon_2$. By the same arguments as above, there exists a neighborhood $\mathcal{U}$ of $(0,0)$ in $\mathbb{R}^2$ such that for each $(a,\alpha) \in (\mathbb{R} \times \mathbb{R}^*) \cap \mathcal{U}$, the distance function to the origin has no critical point on $f^{-1}(\alpha)  \cap \{x_1 =a \} \cap (B_{\epsilon_2} \setminus \mathring{B_{\epsilon_1}})$. Let us choose $(a_3,\alpha_3) \in (\mathbb{R} \times \mathbb{R}^*) \cap \mathcal{U}$ such that $0 < a_3 \le {\rm min}\{a_{\epsilon_1},a_{\epsilon_2} \}$ and $ \alpha_3 \le {\rm min} \{ \alpha_{a_3,\epsilon_1},\alpha_{a_3,\epsilon_2} \}$. Then, by the first case, $f^{-1}(\alpha_3) \cap \{x_1=a_3 \} \cap B_{\epsilon_1}$ is homemorphic to $f^{-1}(\alpha_1) \cap \{x_1=a_1 \} \cap B_{\epsilon_1}$ and $f^{-1}(\alpha_3) \cap \{x_1=a_3 \} \cap B_{\epsilon_2}$ is homemorphic to $f^{-1}(\alpha_2) \cap \{x_1=a_2 \} \cap B_{\epsilon_2}$. But, since the distance function to the origin has no critical point on $f^{-1}(\alpha_3)  \cap \{x_1 =a_3 \} \cap B_{\epsilon_2} \setminus \mathring{B_{\epsilon_1}}$, the fibre
$f^{-1}(\alpha_3) \cap \{x_1=a_3 \} \cap B_{\epsilon_1}$ is homeomorphic to  $f^{-1}(\alpha_3) \cap \{x_1=a_3 \} \cap B_{\epsilon_2}$.
\endproof

Similarly, there exists $\epsilon_0' >0$ such that for $0< \epsilon \le \epsilon_0'$, there exists $b_\epsilon>0$ such that for $0 < a \le b_\epsilon$, the topological type of $f^{-1}(0) \cap \{x_1=a\} \cap B_\epsilon$ and 
$f^{-1}(0) \cap \{x_1=-a\} \cap B_\epsilon$ do not depend on the choice of $(\epsilon,a)$. Therefore we can make the following definition.
\begin{definition}\label{definitiongamma}
{\rm We set 
$$\begin{array}{l}
\gamma^{+,+}= \chi \big( f^{-1}(0) \cap \{x_1=a\} \cap B_\epsilon \big) - \chi \big( f^{-1}(\alpha) \cap \{x_1=a\} \cap B_\epsilon \big) , \cr
\gamma^{+,-}= \chi \big( f^{-1}(0) \cap \{x_1=-a\} \cap B_\epsilon \big) - \chi \big( f^{-1}(\alpha) \cap \{x_1=-a\} \cap B_\epsilon \big) ,\cr
\gamma^{-,+}= \chi \big( f^{-1}(0) \cap \{x_1=a\} \cap B_\epsilon \big) - \chi \big( f^{-1}(-\alpha) \cap \{x_1=a\} \cap B_\epsilon \big) ,\cr
\gamma^{-,-}= \chi \big( f^{-1}(0) \cap \{x_1=-a\} \cap B_\epsilon \big) - \chi \big( f^{-1}(-\alpha) \cap \{x_1=-a\} \cap B_\epsilon \big),
\end{array} $$
where $0< \alpha \ll a \ll \epsilon$. }
\end{definition}
Now we are in position to state the generalization of the Khimshiashvili formula. Remember that $e_1$ satisfies the genericity conditions of Lemmas \ref{Lemma4.1} and \ref{Lemma4.4}.

\begin{theorem}\label{KhimPolar}
Assume that $e_1 \notin \Sigma_1 \cup \Sigma_2$. For $0 < \delta \ll \epsilon \ll 1$, we have 
$$\chi \big( f^{-1}(-\delta) \cap B_\epsilon \big) =1-\lambda^- -\gamma^{-,-} = 1-\lambda^+ -\gamma^{-,+},$$
$$\chi \big( f^{-1}(\delta) \cap B_\epsilon \big) =1-(-1)^n\lambda^- -\gamma^{+,+} = 1-(-1)^n \lambda^+ -\gamma^{+,-}.$$
\end{theorem}
\proof The set of critical points of $x_1$ on $f^{-1}(-\delta) \cap \mathring{B_\epsilon}$ is exactly $\Gamma_{x_1} \cap f^{-1}(-\delta)$. Moreover we know that if $p \in \Gamma_{x_1} \cap f^{-1}(-\delta)$ then $x_1(p) \not= 0$. By Morse theory, we have 
$$\displaylines{ 
\quad \chi \big( f^{-1}(-\delta) \cap B_\epsilon \cap \{x_1 \ge 0 \} \big) - \chi \big( f^{-1}(-\delta) \cap B_\epsilon \cap \{x_1 = 0 \} \big) \hfill \cr
\hfill = \sum_{p \in \Gamma_{x_1} \cap f^{-1}(-\delta) \atop x_1(p)>0 } (-1)^{\mu(p)} , \quad \cr
\quad \chi \big( f^{-1}(-\delta) \cap B_\epsilon \cap \{x_1 \le 0 \} \big) - \chi \big( f^{-1}(-\delta) \cap B_\epsilon \cap \{x_1 = 0 \} \big) \hfill \cr
\hfill = (-1)^{n-1}\sum_{p \in \Gamma_{x_1} \cap f^{-1}(-\delta) \atop x_1(p)<0 } (-1)^{\mu(p)} . \quad \cr
}$$
Here we remark that $f^{-1}(\delta) \cap B_\epsilon$ is a manifold with boundary and $x_1$ may have critical points on the boundary. But  by Lemma 3.7 in \cite{DutertreJofSingProcTrot}, these critical points lie in $\{x_1 \not= 0 \}$ and are outwards-pointing (resp. inwards-pointing) in $\{ x_1 >0 \}$ (resp. $\{x_1 < 0 \}$). That is why they do not appear in the above two formulas. Adding the two equalities and applying the Mayer-Vietoris sequence, we obtain 
$$ \chi \big( f^{-1}(-\delta) \cap B_\epsilon \big) - \chi \big( f^{-1}(-\delta) \cap B_\epsilon \cap \{x_1 = 0 \} \big) =$$
$$\sum_{p \in \Gamma_{x_1} \cap f^{-1}(-\delta) \atop x_1(p)>0 } (-1)^{\mu(p)} + (-1)^{n-1}\sum_{p \in \Gamma_{x_1} \cap f^{-1}(-\delta) \atop x_1(p)<0 } (-1)^{\mu(p)} .$$
Since $\nabla f = f_{x_1} e_1$ on $\Gamma_{x_1}$, it is easy to check that $p$ belongs to $\Gamma_{x_1} \cap \{ f < 0 \} \cap \{x_1 > 0 \}$ if and only if $p$ belongs to $\Gamma_{x_1} \cap \{ f_{x_1} < 0 \} \cap \{x_1 > 0 \}$
and $p$ belongs to $\Gamma_{x_1} \cap \{ f < 0 \} \cap \{x_1 <0 \}$ if and only if $p$ belongs to $\Gamma_{x_1} \cap \{ f_{x_1} > 0 \} \cap \{x_1 < 0 \}$. Let us decompose $\mathcal{B}^+$ into the disjoint union $\mathcal{B}^+= \mathcal{B}^{+,+} \sqcup \mathcal{B}^{+,-}$ where $\mathcal{B}^{+,+}$ (resp. $\mathcal{B}^{+,-}$) is the set of half-branches of $\mathcal{B}^+$ on which $x_1 >0$ (resp. $x_1 < 0$). Similarly we can write  $\mathcal{B}^-= \mathcal{B}^{-,+} \sqcup \mathcal{B}^{-,-}$. Combining Lemma \ref{Lemma4.6} and Lemma \ref{Lemma4.7}, we can rewrite the above equality in the following form:
$$\chi \big( f^{-1}(-\delta) \cap B_\epsilon \big) - \chi \big( f^{-1}(-\delta) \cap B_\epsilon \cap \{x_1 = 0 \} \big) = -\sum_{{\rm \bf b} \in \mathcal{B}^{-,+} } \sigma({\rm \bf b}) - \sum_{{\rm \bf b} \in \mathcal{B}^{+,-} } \sigma({\rm \bf b}). \eqno(1)$$
Since $(-\delta,0)$ is a regular value of $(f,x_1)$ then there exists $a_\delta >0$ such that for $0< a \le a_\delta$, $(-\delta,\pm a)$ are regular value of $(f,x_1)$ and 
$$\chi \big( f^{-1}(-\delta) \cap \{ x_1 = \pm a \} \cap B_\epsilon \big) = \chi \big( f^{-1}(-\delta) \cap \{ x_1 = 0 \} \cap B_\epsilon \big).$$ 

Let us fix $a$ such that $0< a \le a_\delta$ and let us relate $\chi \big( f^{-1}(-\delta) \cap \{ x_1 = - a \} \cap B_\epsilon \big)$ to $\chi \big( f^{-1}(\alpha) \cap \{ x_1 =  -a \} \cap B_\epsilon \big)$ where $0 < \alpha \ll a$. 
Note that the set of critical points of $f$ on $\{x_1 = -a \} \cap \mathring{B_\epsilon}$ is exactly $\Gamma_{x_1} \cap \{ x_1 = -a \}$. Moreover this set of critical points is included in $\{ f > - \delta \}$. Indeed, if it is not the case, then there is a half-branch of $ \Gamma_{x_1}$ that intersects $\{x_1 = -a\}$ on $\{ f \le -\delta \}$. But since $x_1$ and $f$ are negative on this branch, this would imply that $\Gamma_{x_1}$ intersects $\{ f =- \delta \}$ on $\{-a \le x_1 < 0 \}$, which is not possible for $a \le a_\delta$. 

Now let us look at the critical points of $f$ on $\{x_1=-a\} \cap S_\epsilon$. In the proof of Lemma \ref{Lemma4.10}, we established the existence of a neighborhood $\mathcal{U}_\epsilon$ of $(0,0)$ in $\mathbb{R}^2$ such that for each $(a,\alpha) \in (\mathbb{R} \times \mathbb{R}^*) \cap \mathcal{U}_\epsilon$, the fibre $f^{-1}(\alpha) \cap \{x_1 = a \}$ intersects $S_\epsilon$ transversally. Therefore we can choose $\delta$ such that the critical points of $f$ on $\{x_1 = 0 \} \cap S_\epsilon \cap \{ f \not= 0 \}$ lie in $\{ \vert f \vert > \delta \}$. Moreover by a Curve Selection Lemma argument, they are outwards-pointing in $\{ f > \delta \}$ and inwards-pointing in $\{ f < -\delta \}$. So, if $a$ is small enough, then the critical points of $f$ on $\{x_1 = -a \} \cap S_\epsilon \cap \{ f \not= 0 \}$, lying in $\{ \vert f \vert > \delta \}$, are outwards-pointing (resp. inwards-pointing) in $\{ f  > \delta \}$ (resp. $\{ f < -\delta \}$). By Morse theory, we find that
$$\displaylines{ \quad \chi \big( \{ f \le -\delta \} \cap \{x_1 = -a \} \cap B_\epsilon \big) -\chi \big( \{ f = -\delta \} \cap \{x_1 = -a \} \cap B_\epsilon \big) = 0, \hfill \cr }$$
$$\displaylines{
\quad \chi \big( \{-\delta \le  f \le -\alpha \} \cap \{x_1 = -a \} \cap B_\epsilon \big) -\chi \big( \{ f = -\delta \} \cap \{x_1 = -a \} \cap B_\epsilon \big) \hfill \cr
\hfill = \sum_{p  \in \Gamma_{x_1} \cap \{x_1 =-a \} \atop f(p)<0} (-1)^{\theta(p)}, \quad \cr 
\quad \chi \big( \{ f \ge \alpha \} \cap \{x_1 = -a \} \cap B_\epsilon \big) -\chi \big( \{ f = \alpha \} \cap \{x_1 = -a \} \cap B_\epsilon \big) \hfill \cr
\hfill = \sum_{p  \in \Gamma_{x_1} \cap \{x_1 =-a \} \atop f(p)>0} (-1)^{\theta(p)}. \quad \cr
}$$
By the Mayer-Vietoris sequence, we have that 
$$\displaylines{
\quad 1 = \chi \big( \{x_1=-a \} \cap B_\epsilon \big) = \chi \big( \{ f \le -\delta \} \cap \{x_1 = -a \} \cap B_\epsilon \big) \hfill \cr
\quad \quad + \chi \big( \{-\delta \le  f \le -\alpha \} \cap \{x_1 = -a \} \cap B_\epsilon \big)
- \chi \big( \{ f = -\delta\} \cap \{x_1 = -a \} \cap B_\epsilon \big) \hfill \cr
\quad \quad \quad - \chi \big( \{f = -\alpha \} \cap \{x_1 = -a \} \cap B_\epsilon \big) + \chi \big( \{-\alpha \le  f \le \alpha \} \cap \{x_1 = -a \} \cap B_\epsilon \big) \hfill \cr
\quad \quad \quad \quad + \chi \big( \{ f \ge \alpha \} \cap \{x_1 = -a \} \cap B_\epsilon \big) - \chi \big( \{ f =\alpha \} \cap \{x_1 = -a \} \cap B_\epsilon \big). \hfill \cr
}$$
Using the fact that the inclusion $$\{f=0\} \cap \{x_1 = -a \} \cap B_\epsilon \subset \{-\alpha \le  f \le \alpha \} \cap \{x_1 = -a \} \cap B_\epsilon$$ is a homotopy equivalence and applying the above equalities, we get 
$$1=  \sum_{p  \in \Gamma_{x_1} \cap \{x_1 =-a \} \atop f(p)<0} (-1)^{\theta(p)} +  \sum_{p  \in \Gamma_{x_1} \cap \{x_1 =-a \} \atop f(p)>0} (-1)^{\theta(p)} + \chi \big( f^{-1}(-\delta) \cap \{x_1=-a \} \cap B_\epsilon \big) $$
$$+ \chi \big( f^{-1}(0) \cap \{x_1=-a \} \cap B_\epsilon \big) - \chi \big( f^{-1}(-\alpha) \cap \{x_1=-a \} \cap B_\epsilon \big).$$ 
By Lemma \ref{Lemma4.7}, we can rewrite this equality in the following form:
$$1= -\sum_{{\rm \bf b} \in \mathcal{B}^{+,-}} \sigma ({\rm \bf b}) + \sum_{{\rm \bf b} \in \mathcal{B}^{-,-}} \sigma ({\rm \bf b}) + \chi \big( f^{-1}(-\delta) \cap \{x_1=-a\} \cap B_\epsilon \big) + \gamma^{-,-} . \eqno(2)$$
Combining $(1)$ and $(2)$, we obtain the first equality of the statement. The second one is obtained replacing $-a$ with $a$ in the above discussion. The third and fourth ones are obtained replacing $f$ with $-f$. 
\endproof

\begin{remark}
{\rm \begin{enumerate}
\item If $f$ has an isolated critical point at the origin then we recover the Khimshiashvili formula because $$\gamma^{-,-}=\gamma^{+,+}=\gamma^{+,-}=\gamma^{-,+}=0$$ and $\lambda^+=\lambda^-= {\rm deg}_0 \nabla f$.
\item If we denote by $S_0$ the stratum that contains $0$ and if we assume that ${\rm dim} \ S_0 >0$, then by the Thom $(a_f)$-condition, the polar curve $\Gamma_v$ is empty in a neighborhood of $0$ if $v \notin S^{n-1} \cap (T_x S_0)^\perp$. Then applying Equality (1) of the previous proof, we recover Corollary \ref{Corollary3.5}. Actually, we can say more about the relation between the topologies of $f^{-1}(\pm \delta) \cap B_\epsilon$ and $f^{-1}(\pm \delta) \cap B_\epsilon \cap \{ v^* = 0 \}$. As mentionned in the proof of Theorem \ref{KhimPolar}, the critical points of $v^*$ restricted to $f^{-1}(\pm \delta) \cap S_\epsilon$ lie in $\{v^* \not= 0 \}$ and are outwards-pointing (resp. inwards-pointing) in $\{v^* > 0 \}$ (resp. $\{ v^* < 0 \}$). So we can apply the arguments of the proof of Theorem 6.3 in \cite{DutertreAraujoIMRN} to get that $f^{-1}(\pm \delta) \cap B_\epsilon$ is homeomorphic to $f^{-1}(\pm \delta) \cap B_\epsilon \cap \{ v^* = 0 \} \times [-1,1]$.
\end{enumerate}}
\end{remark}

\section{One dimensional critical locus and a real L\^e-Iomdine formula}
In this section, we apply the results of Section 4 to the case of a one-dimensional singular set, in order to establish a real version of the L\^e-Iomdine formula. 

Let $f : (\mathbb{R}^n,0) \to (\mathbb{R},0)$ be a definable function-germ of class $C^r$, $2 \le r$. We assume that ${\rm dim}\ \Sigma_f =1$. 
In the neighborhood of the origin, the partition 
$$ \big( f^{-1}(0) \setminus \Sigma_f, \Sigma_f \setminus \{0\}, \{0\} \big)$$
gives a Whitney stratification of $f^{-1}(0)$ which satisfies the Thom $(a_f)$-condition, because the points where the Whitney conditions and the Thom $(a_f)$-condition may fail form a $0$-dimensional definable set of $\Sigma_f \setminus \{0\}$. 
Let $\mathcal{C}$ be the set of half-branches of $\Sigma_f$, i.e., the set of connected components of $\Sigma_f \setminus \{0\}$. 
\begin{lemma}\label{Lemma5.1}
There exists a definable set $\Sigma_3 \subset S^{n-1}$ of positive codimension such that if $v \notin \Sigma_3$, $v^*$ does not vanish on $\Sigma_f \setminus \{0\}$ in a neighborhood of the origin.
\end{lemma}
\proof Let $c \in \mathcal{C}$. If $v^*$ vanishes on ${\rm \bf c}$ in a neighborhood of the origin then, if $u \not= 0 $ is on $C_0 {\rm \bf c}$ (the tangent cone at ${\rm \bf c}$ at the origin) then $v^*(u)=0$ and so $v \in u^\perp$. So if $v \notin \cup_{{\rm \bf c} \in \mathcal{C}} (C_0 {\rm \bf c})^\perp$ then $v^*$ does not vanish on $\Sigma_f \setminus \{0\}$. But $(\cup_{{\rm \bf c} \in \mathcal{C}} (C_0 {\rm \bf c})^\perp) \cap S^{n-1}$ has dimension less than or equal to $n-2$. \endproof

From now on, we assume that $e_1 \in S^{n-1}$ is generic, i.e., $e_1 \notin \Sigma_1 \cup \Sigma_2 \cup \Sigma_3$. Since $e_1 \notin \Sigma_3$, there exists $a_1 >0$ such that for $0 < a \le a_1$, $x_1^{-1}(\pm a)$ intersects $\Sigma_f \setminus \{0\}$ transversally and so, the points in $x_1^{-1}(\pm a) \cap (\Sigma_f \setminus \{0\})$ are isolated critical points of $f_{\vert \{x_1 = \pm a \}}$. For $q \in x_1^{-1}(\pm a) \cap (\Sigma_f \setminus \{0\})$, we denote by ${\rm deg}_q \nabla f_{\vert x_1^{-1}(\pm a)}$ the topological degree of the mapping $\frac{\nabla f_{\vert x_1^{-1}(\pm a)}}{\vert \nabla f_{\vert x_1^{-1}(\pm a) }\vert} : x_1^{-1}(\pm a) \cap S_{\epsilon'} (q) \to S^{n-2}$, where $S_{\epsilon'} (q)$ is the sphere centered at $q$ of radius $\epsilon'$ with $0< \epsilon' \ll 1$.

Let us write $\mathcal{C}=\mathcal{C}^+ \sqcup \mathcal{C}^-$ where $\mathcal{C}^+$ (resp. $\mathcal{C}^-$) is the set of half-branches of $\mathcal{C}$ on which $x_1 >0$ (resp. $x_1 <0$).
\begin{lemma}
Let ${\rm \bf c} \in \mathcal{C}^+$. There exists $a_{\rm \bf c} >0$ such that the function $a \mapsto {\rm deg}_q \nabla f_{\vert x_1^{-1}(a)}$, where $\{q\} = {\rm \bf c} \cap \{x_1 = a\}$, is constant on $]0,a_{\rm \bf c}]$.
\end{lemma}
\proof It is enough to prove that there exists an interval  $]0,a_{\rm \bf c}]$ on which the function $a \mapsto {\rm deg}_q \nabla f_{\vert x_1^{-1}(a)}$ is locally constant. Let $d : \mathbb{R}^n \to \mathbb{R}$ be the distance function to ${\rm \bf c}$. It is a continuous definable function and there exists an open definable neighbourhood  $\mathcal{U}$ of ${\rm \bf c}$ such that $d$ is smooth on $\mathcal{U} \setminus {\rm \bf c}$. Moreover we can assume that $d$ is a (stratified) submersion on $\{f \le 0 \} \cap (\mathcal{U} \setminus {\rm \bf c})$. 

Let $\pi : \{f \le 0 \} \cap (\mathcal{U} \setminus {\rm \bf c}) \to \mathbb{R}^2$ be the mapping defined by $\pi(p)=(x_1(p),d(p))$ and let $\Delta \subset \mathbb{R}^2$ be its (stratified) discriminant. It is a definable curve 
included in $\mathbb{R} \times \mathbb{R}^*$ and so $\bar{\Delta} \cap (\mathbb{R} \times \{0\})$ is a finite number of points. Let us choose $a_{\rm \bf c} > 0$ such that 
$$ a_{\rm \bf c} < {\rm min} \left\{ x_1 (u) \ \vert \ u \in \bar{\Delta} \cap (\mathbb{R}^* \times \{0\} ) \right\}.$$
If $0 < a \le a_{\rm \bf c}$, then there exists $t>0$ and $\epsilon > 0$ such that  $]a-t,a+t [ \times ]0,\epsilon[$ does not meet $\Delta$. Hence the function 
$$\begin{array}{ccc}
]a-t,a+t [ \times ]0,\epsilon[  & \rightarrow & \mathbb{R} \cr
 (a',\epsilon') & \mapsto & \chi \big( \{ f \le 0 \} \cap \{x_1 = a' \} \cap \{ d =\epsilon ' \} \big) \cr
\end{array}$$
is constant.  Therefore by Corollary \ref{ArnoldWall}, the function $a ' \mapsto  {\rm deg}_{q'} \nabla f_{\vert x_1^{-1}(a')}$ is constant on $]a-t,a+t[$. \endproof
Of course, a similar result is valid for ${\rm \bf c} \in \mathcal{C}^-$. If ${\rm \bf c} \in \mathcal{C}$, let us denote by $\tau({\rm \bf c})$ the value that the function $a \mapsto {\rm deg}_q \nabla f_{\vert x_1^{-1}(a)}$, $\{q\} = {\rm \bf c} \cap \{x_1 = a\}$, takes close to the origin.
\begin{definition}
{\rm We set $\gamma^+ = \sum_{{\rm \bf c} \in \mathcal{C}^+} \tau ({\rm \bf c})$ and $ \gamma^- = \sum_{{\rm \bf c} \in \mathcal{C}^-} \tau ({\rm \bf c})$.}
\end{definition}
In this setting, Theorem \ref{KhimPolar} admits the following formulation.
\begin{theorem}\label{KhimOneDim}
Assume that $e_1 \notin \Sigma_1 \cup \Sigma_2 \cup \Sigma_3$. For $0 < \delta \ll \epsilon \ll 1$, we have
$$\chi \big( f^{-1}(-\delta) \cap B_\epsilon \big) = 1-\lambda^--\gamma^-=1-\lambda^++\gamma^+,$$
$$\chi \big( f^{-1}(\delta) \cap B_\epsilon \big) = 1-(-1)^n (\lambda^+-\gamma^-)=1-(-1)^n (\lambda^--\gamma^+).$$
\end{theorem}
\proof Since $\Gamma_{x_1} \cap  f^{-1}(0) = \emptyset$, the critical points of $f_{\vert \{x_1 =\pm a \} }$ in $f^{-1}([-\alpha,\alpha])$, $0 < \alpha \ll \delta \ll \epsilon$, are exactly the points in $\Sigma_f \cap \{x_1 = \pm a \}$. An easy adaptation of the proof of the Khimshiashvili formula (Theorem \ref{KhimshiashviliFormula}) gives that 
$$\gamma^{-,-}= \gamma^- ,\ \gamma^{-,+}=\gamma^+,\ 
\gamma^{+,+}= (-1)^{n-1} \gamma^+ \hbox{ and } \gamma^{+,-}=(-1)^{n-1} \gamma^-.$$
\endproof
We remark that $\lambda^-+\gamma^-=\lambda^++\gamma^+$. Moreover, if $n$ is even, we have 
$$\chi \big( f^{-1}(\delta) \cap B_\epsilon \big) - \chi \big( f^{-1}(-\delta) \cap B_\epsilon \big) = \gamma^+ + \gamma^-,$$
$$\chi \big( f^{-1}(\delta) \cap B_\epsilon \big) + \chi \big( f^{-1}(-\delta) \cap B_\epsilon \big) = 2 -(\lambda^+ + \lambda^-),$$
and if $n$ is odd, we have
$$\chi \big( f^{-1}(\delta) \cap B_\epsilon \big) - \chi \big( f^{-1}(-\delta) \cap B_\epsilon \big) = \lambda^+ + \lambda^-,$$
$$\chi \big( f^{-1}(\delta) \cap B_\epsilon \big) + \chi \big( f^{-1}(-\delta) \cap B_\epsilon \big) = 2 -(\gamma^+ + \gamma^-). $$
Therefore the two sums $\sum_{{\rm \bf b} \in \mathcal{B}} \sigma({\rm \bf b})$ and $\sum_{{\rm \bf c} \in \mathcal{C}} \tau({\rm \bf c})$ do not depend on the generic choice of linear function that we used to define them. Moreover, applying Lemma \ref{MilnorfibreLink}, we get that if $n$ is even, $\chi ({\rm Lk}(\{ f=0 \}))= 2 -(\lambda^+ + \lambda^-)$ and if $n$ is odd,  $\chi ({\rm Lk}(\{ f=0 \}))=\gamma^+ + \gamma^-$. 

Let us give an example. Let $f(x,y,z)=y^2-zx^b$, $b >1$ (see \cite{MasseyLeCycles}, Example 2.2). This polynomial is weighted-homogeneous but we cannot apply Corollary \ref{FibreMilnorHomogeneous}, for $b$ may be arbitrary large. Then $\Sigma_f = \{ (0,0,z) \ \vert \ z \in \mathbb{R} \}$. 

Let $v=(1,1,1)$ so that $v^*(x,y,z)=x+y+z$. We have to check that $v$ satisfies the conclusions of Lemma \ref{Lemma4.1}, Lemma \ref{Lemma4.4} and Corollary \ref{Corollary4.5}, and Lemma \ref{Lemma5.1}. A straightforward computation shows that 
$$\Gamma_v = \left\{ (x,-\frac{x^b}{2},\frac{x}{b})\ \vert \ x \not= 0 \right\}.$$
Since $\Gamma_v \cap \{ v^* = 0 \} = \emptyset$, we see that $\{v^*=0\}$ intersects the stratum $f^{-1}(0) \setminus \Sigma_f$ transversally. Moreover, since $v^*$ does not vanish on $\Sigma_f \setminus \{0\}$, $\{v^*=0 \}$ intersects the stratum $\Sigma_f \setminus \{0\}$ transversally and so $v$ satisfies the conclusion of Lemma \ref{Lemma4.1} (and of Lemma \ref{Lemma5.1} as well). It is clear that $\Gamma_v$ is a curve in the neighborhood of the origin. In order to check that the conclusion of Corollary \ref{Corollary4.5} holds, thanks to the computations of Lemmas \ref{Lemma4.6} and \ref{Lemma4.7}, it is enough to check that ${\rm det}[\nabla f_x,\nabla f_y,\nabla f_z ]$ does not vanish on $\Gamma_v$. But 
$${\rm det}[\nabla f_x,\nabla f_y,\nabla f_z ](x,y,z)= -2 b^2 x^{2b-2},$$
and so the conclusion of Corollary \ref{Corollary4.5} holds. Moreover, since 
$$\frac{\partial f}{\partial v}(x,y,z)= -bx^{b-1}z+2y-x^b,$$
we easily compute that $\lambda^+=\lambda^-=-1$ if $b$ is odd and that $\lambda^+=0$ and $\lambda^- = -2$ if $b$ is even. 

It remains to compute $\gamma^+$ and $\gamma^-$. But $\gamma^+$ is the local topological degree at $(0,0)$ of the function $f(x,y,a-x-y)$, $a>0$, that is the local topological degree at $(0,0)$ of the function 
$$(x,y) \mapsto y^2-ax^b + x^{b+1} + yx^b.$$
Then it is not difficult to see that $\gamma^+=-1$ if $b$ is even and $\gamma^+= 0$ if $b$ is odd. Similarly $\gamma^-=1$ if $b$ is even and $\gamma^-=0$ if $b$ is odd. Therefore, applying Theorem \ref{KhimOneDim} and Lemma \ref{MilnorfibreLink}, we obtain that 
$$\chi \big( f^{-1}(-1) \big) = 2 , \ \chi \big( f^{-1}(1) \big) = 0 
\hbox{ and } 
\chi ({\rm Lk}(\{ f=0 \})) = 0.$$ 

In the rest of the section, we will apply Theorem \ref{KhimOneDim} to establish a real version of the L\^e-Iomdine formula. From now on, we assume that the structure is polynomially bounded. 
\begin{lemma}
There exists $n_0 \in \mathbb{N}$ such that 
$$\vert f_{x_1}(p) \vert > \vert x_1 (p) \vert^{n_0} \hbox{ and } \vert f(p) \vert > \vert x_1 (p) \vert^{n_0},$$
for $p \in \Gamma_{x_1}$ close to the origin.
\end{lemma}
\proof For $u>0$ small, we define $\beta(u)$ by 
$$\beta (u)= {\rm inf} \big\{ \vert f_{x_1}(p) \vert \ \vert \ p \in \Gamma_{x_1} \cap \{ \vert x_1 \vert = u \} \big\}.$$
It is well defined because $\Gamma_{x_1} \cap \{ \vert x_1 \vert = u \} $ is finite and $\beta (u) >0$. The function $\beta$ is definable and so is the function $\alpha (R)=\beta (\frac{1}{R})$, defined for $R>0$ sufficiently big. Then there exists $n_0 \in \mathbb{N}$ such that $\frac{1}{\alpha(R)} < R^{n_0}$ for $R>0$ sufficiently big. This implies that $\beta (\frac{1}{R}) > \frac{1}{R^{n_0}}$, i.e., $\beta (u) > u^{n_0}$ for $u>0$ sufficiently small. Hence for $p \in \Gamma_{x_1}$ sufficiently close to the origin, we have
$$\vert f_{x_1}(p) \vert > \vert x_1(p) \vert^{n_0}.$$
A similar proof works for the second equality because $f$ and $x_1$ do not vanish on $\Gamma_{x_1}$. \endproof

Let us fix $k \in \mathbb{N}$ with $k > n_0 +1$ and let us set $g (x)=f(x)+ x_1^k$.
\begin{lemma}
The function $g$ has an isolated critical point at the origin.
\end{lemma}
\proof A point $p$ belongs to $(\nabla g)^{-1}(0)$ if and only if 
$$\frac{\partial f}{\partial x_1}(p) + k x_1^{k-1}(p)=0 \hbox{ and }  \frac{\partial f}{\partial x_i}(p) =0 \hbox{ for } i \ge 2.$$
Let us suppose first that $p \in \Sigma_f \setminus \{ 0 \}$. This implies that $x_1 (p)=0$. Since $x_1$ does not vanish on $\Sigma_f \setminus \{0\}$ close to the origin, this case is not possible. Let us suppose now that $p \notin \Sigma_f$. Then $p$ belongs to $\Gamma_{x_1}$ and so $x_1(p) \not= 0$ and $f_{x_1}(p) \not= 0$. By the previous lemma, $\vert f_{x_1}(p) \vert > \vert x_1(p) \vert^{n_0}$ which implies that $k \vert x_1(p) \vert^{k-1} > \vert x_1 (p) \vert^{n_0}$, and so $\vert x_1 (p) \vert^{k-n_0-1} > \frac{1}{k}$ in the neighborhood of the origin. This is impossible by the choice of $k$. The only possible case is when $p$ is the origin. \endproof

The previous lemma unables us to use the Khimshiashvili formula to compute the Euler characteristic of the Milnor fibre of $g$.  We will relate ${\rm deg}_0  \nabla g$ to the indices $\lambda^+$, $\lambda^-$, $\gamma^+$ and $\gamma^-$. Before that we need some auxiliary results. Let 
$$\Gamma_{x_1}(g) =  \left\{ x \in \mathbb{R}^n \setminus \Sigma_g \ \vert \ {\rm rank}(\nabla g (x), e_1 ) < 2 \right\}.$$
\begin{lemma}
We have $\Gamma_{x_1}(g) \cap \{g =0 \} =\emptyset$.
\end{lemma}
\proof If it is not the case this implies that the following set 
$$\{g=0\} \cap \{ g_{x_2}=\ldots=g_{x_n}=0 \} \cap \{x_1 =0\} \setminus \{0\}$$
is not empty in the neighbourhood of the origin. Therefore the set 
$$\{f=0\} \cap \{ f_{x_2}=\ldots=f_{x_n}=0 \} \cap \{x_1 =0\} \setminus \{0\}$$
is not empty in the neighbourhood of the origin. But this is not possible because $\{f=0\} \cap \{ f_{x_2}=\ldots=f_{x_n} \}= \Sigma_f$ and $\{x_1 = 0 \} \setminus \{0\} \cap \Sigma_f = \emptyset$. \endproof
\begin{lemma}
The set $\Gamma_{x_1}(g)$ admits the following decomposition:
$$ \Gamma_{x_1}(g) = \Gamma_{x_1} \sqcup (\Sigma_f \setminus \{0\}).$$
\end{lemma}
\proof We see that $p \in \Gamma_{x_1}(g)$ if and only if $g_{x_2}(p) =\ldots= g_{x_n}(p)=0$ and $g_{x_1} (p) \not= 0$. Since $g_{x_i}(p)=f_{x_i}(p)$, $i=2,\ldots,n$, it is clear that $\Gamma_{x_1} \sqcup (\Sigma_f \setminus \{0\}) \subset \Gamma_{x_1}(g)$. If $p \in \Gamma_{x_1}(g)$ then $f_{x_2}(p)=\ldots=f_{x_n}(p)=0$ and $f_{x_1}(p) + k x_1(p)^{k-1} \not= 0$. If $f_{x_1}(p) \not= 0$ then $p \in \Gamma_{x_1}$. If $f_{x_1}(p)=0$ then $p \in \Sigma_f \setminus \{0\}$. \endproof

\begin{lemma}\label{Lemma5.9}
If $p \in \Gamma_{x_1}$, then ${\rm det} \big[ \nabla g_{x_1}(p),\ldots, \nabla g_{x_n}(p) \big] \not= 0$ and 
$$\hbox{\rm sign det}  \big[ \nabla g_{x_1}(p),\ldots, \nabla g_{x_n}(p) \big] =
\hbox{\rm sign det}  \big[ \nabla f_{x_1}(p),\ldots, \nabla f_{x_n}(p) \big] .$$
\end{lemma}
\proof Let $p \in \Gamma_{x_1}$. We have $g_{x_1}(p)= f_{x_1}(p) + kx_1^{k-1}$. By the choice of $k$, this implies that ${\rm sign}\  g_{x_1}(p)={\rm sign}\ f_{x_1}(p)$. Using the computations of Lemmas \ref{Lemma4.6} and \ref{Lemma4.7}, we see that
$$ \hbox{sign det}  \big[ \nabla f_{x_1}(p),\ldots, \nabla f_{x_n}(p) \big] = {\rm sign} \left( x_1(p) f_{x_1}(p)  \frac{\partial (x_1,f_{x_2},\ldots,f_{x_n})}{\partial (x_1,\ldots,x_n)} (p) \right).$$
But $$\frac{\partial (x_1,f_{x_2},\ldots,f_{x_n})}{\partial (x_1,f_{x_2},\ldots,x_n)} (p) =\frac{\partial (x_1,g_{x_2}\ldots,g_{x_n})}{\partial (x_1,\ldots,x_n)} (p)$$ so $ \frac{\partial (x_1,g_{x_2},\ldots,g_{x_n})}{\partial (x_1,\ldots,x_n)} (p) \not= 0$. Since $x_1(p) g_{x_1}(p) \not= 0$, we obtain that $${\rm det} \big[ \nabla g_{x_1}(p),\ldots, \nabla g_{x_n}(p) \big] \not= 0$$ and since ${\rm sign}\ g_{x_1}(p)= {\rm sign}\ f_{x_1}(p)$, we conclude that 
$$\hbox{sign det}  \big[ \nabla g_{x_1}(p),\ldots, \nabla g_{x_n}(p) \big] =
\hbox{sign det}  \big[ \nabla f_{x_1}(p),\ldots, \nabla f_{x_n}(p) \big] .$$  \endproof

\begin{lemma}\label{Lemma5.10} 
Assume that $k$ is even. If $q \in \Sigma_{f} \setminus \{0\}$ is close enough to the origin and $x_1(q)=a$, then ${\rm deg}_q \nabla f _{\vert \{x_1 = a \}}$ is equal to ${\rm deg}_q (\nabla g- \nabla g(q))$, where ${\rm deg}_q (\nabla g- \nabla g(q))$ is the topological degree of the mapping $\frac{\nabla g- \nabla g(q)}{\vert \nabla g- \nabla g(q) \vert} : S_{\epsilon'} (q) \to S^{n-1}$ with $0 < \epsilon' \ll 1$.
\end{lemma}
\proof We have that  ${\rm deg}_q \nabla f _{\vert \{x_1 = a \}}$ is equal to the topological degree of the mapping $\frac{W}{\vert W \vert} : S_{\epsilon'}(q) \to S^{n-1}$ where $W=(x_1-a,f_{x_2},\ldots,f_{x_n})$. But 
$$\nabla g- \nabla g(q) =(f_{x_1}+k x_1^{k-1}-ka^{k-1},f_{x_2},\ldots,f_{x_n} )$$ and so, since $f_{x_1}(q)=0$ and $k-1$ is odd, there exists a small neighborhood of $q$ on which $f_{x_1}+k x_1^{k-1}-ka^{k-1}$ and $x_1-a$ have the same sign. If $\epsilon'$ is small enough, then the mappings $\frac{\nabla g- \nabla g(q)}{\vert \nabla g- \nabla g(q) \vert}$ and $\frac{W}{\vert W \vert}$ are homotopic on $S_{\epsilon'} (q)$. Hence the two topological degrees are equal. \endproof

\begin{proposition}
If $k$ is odd then ${\rm deg}_0 \nabla g =\lambda^-=\lambda^+ +\gamma^+ - \gamma^-$. 
If $k$ is even then ${\rm deg}_0 \nabla g =\lambda^-+  \gamma^-=\lambda^+ +\gamma^+$. 
\end{proposition}
\proof Let $\eta > 0$ be a small real number. The set $(\nabla g)^{-1}(-\eta,0,\ldots,0)$ is finite because $\Gamma_{x_1}(g)$ is one-dimensional. Let us write $$(\nabla g)^{-1}(-\eta,0,\ldots,0) = \{p_1,\ldots,p_s\} \cup \{q_1,\ldots,q_r \},$$ where $$ \{p_1,\ldots,p_s\} = (\nabla g)^{-1}(-\eta,0,\ldots,0) \cap \Gamma_{x_1}$$ and 
$$\{q_1,\ldots,q_r \}=  (\nabla g)^{-1}(-\eta,0,\ldots,0) \cap  (\Sigma_f \setminus \{0\}).$$ Therefore we have 
$${\rm deg}_0 \nabla g =\sum_{i=1}^s {\rm deg}_{p_i} (\nabla g- \nabla g(p_i)) + \sum_{j=1}^r {\rm deg}_{q_j} (\nabla g- \nabla g(q_j)).$$
If $k$ is odd, $  (\nabla g)^{-1}(-\eta,0,\ldots,0) \cap  (\Sigma_f \setminus \{0\})$ is empty and, by the choice of $k$, 
$$ (\nabla g)^{-1}(-\eta,0,\ldots,0) \cap \Gamma_{x_1} =  (\nabla g)^{-1}(-\eta,0,\ldots,0) \cap \left[  \cup_{b \in \mathcal{B}^-} b \right] .$$
Using Lemma \ref{Lemma5.9}, we conclude that 
$${\rm deg}_0 \nabla g = \lambda^- = \lambda^++\gamma^+ -\gamma^-.$$
If $k$ is even then 
$$ (\nabla g)^{-1}(-\eta,0,\ldots,0) \cap ( \Sigma_f \setminus \{0\}) =  (\nabla g)^{-1}(-\eta,0,\ldots,0) \cap \left[  \cup_{c \in \mathcal{C}^-} c \right] .$$
Using Lemma \ref{Lemma5.10}, we conclude that ${\rm deg}_0 \nabla g = \lambda^- + \gamma^-= \lambda^+ + \gamma^+$.
\endproof
Now we are in position to formulate the real version of the L\^e-Iomdine formula.
\begin{theorem}\label{RealLeIomdine}
Assume that $e_1 \notin \Sigma_1 \cup  \Sigma_2 \cup \Sigma_3$ and that $k > n_0 +1$. For $0< \delta \ll \epsilon \ll 1$, we have 
\begin{itemize}
\item[-] if $k$ is odd, 
$$\chi \big( g^{-1}(-\delta) \cap B_\epsilon \big) =\chi \big( f^{-1}(-\delta) \cap B_\epsilon \big) + \gamma^-,$$
$$\chi \big( g^{-1}(\delta) \cap B_\epsilon \big) =\chi \big( f^{-1}(\delta) \cap B_\epsilon \big) +(-1)^{n-1} \gamma^+,$$
\item[-] if $k$ is even, 
$$\chi \big( g^{-1}(-\delta) \cap B_\epsilon \big) =\chi \big( f^{-1}(-\delta) \cap B_\epsilon \big) ,$$
$$\chi \big( g^{-1}(\delta) \cap B_\epsilon \big) =\chi \big( f^{-1}(\delta) \cap B_\epsilon \big) +(-1)^{n-1} (\gamma^+ + \gamma^-).$$
\end{itemize}
\end{theorem}
\proof We know that $\chi \big( g^{-1}(-\delta) \cap B_\epsilon \big)=1-{\rm deg}_0 \nabla g$. If $k$ is odd, this gives 
$$\chi \big( g^{-1}(-\delta) \cap B_\epsilon \big)= 1 -\lambda^-= \chi \big( f^{-1}(-\delta) \cap B_\epsilon \big) + \gamma^-.$$
If $k$ is even, this gives 
$$\chi \big( g^{-1}(-\delta) \cap B_\epsilon \big)= 1 -\lambda^--\gamma^-= \chi \big( f^{-1}(-\delta) \cap B_\epsilon \big) .$$
We know that $\chi \big( g^{-1}(\delta) \cap B_\epsilon \big)=1- (-1)^n {\rm deg}_0 \nabla g$. So if $k$ is odd, then
$$\chi \big( g^{-1}(\delta) \cap B_\epsilon \big)= 1 -(-1)^n \lambda^-= \chi \big( f^{-1}(\delta) \cap B_\epsilon \big) -(-1)^n \gamma^+.$$
If $k$ is even, we get that 
$$\displaylines{
\quad	\chi \big( g^{-1}(\delta) \cap B_\epsilon \big)= 1 -(-1)^n( \lambda^- + \gamma^-) \hfill \cr
\qquad \qquad \qquad \qquad \qquad = 1- (-1)^n (\lambda^- - \gamma^+) - (-1)^n (\gamma^+ + \gamma^-) \hfill \cr
\hfill   =  \chi \big( f^{-1}(\delta) \cap B_\epsilon \big) -(-1)^n( \gamma^+ + \gamma^-). \quad \quad \quad  \cr
}$$
\endproof

\end{document}